\documentclass[a4paper,english,12pt]{article} 
\usepackage{overpic}
\usepackage{todonotes}
\usepackage{mystyle}
\usepackage{mymath} %
\usepackage{mylisting} %
\usepackage{authblk} 
\usepackage{blindtext} 
\usepackage{multicol}
\usepackage{sectsty} 
\usepackage{caption}
\usepackage{geometry}
\usepackage{subcaption}

\crefname{lstlisting}{Code}{listings}
\Crefname{lstlisting}{Code}{Listings}
\usepackage[mathlines]{lineno} 
\usepackage{doi}
\usepackage[section]{placeins}
\usepackage[normalem]{ulem}

\newtheorem{theorem}{Theorem}
\newcommand*{\TitleFont}{%
      \usefont{\encodingdefault}{\rmdefault}{b}{}%
      \fontsize{16}{15}%
      \selectfont} 
\sectionfont{\fontsize{14}{2}\selectfont}
\subsectionfont{\fontsize{12}{2}\selectfont} 

\title{\TitleFont Automatic shape derivatives for transient PDEs in FEniCS and Firedrake}
\author[1,2,*]{Jørgen S. Dokken}
\author[2]{Sebastian K. Mitusch}
\author[2]{Simon W. Funke}
\affil[1]{Department of Engineering, University of Cambridge, CB2 1PZ,
United Kingdom}
\affil[2]{Simula Research Laboratory, 1325 Lysaker, Norway}
\affil[*]{Email: jsd55@cam.ac.uk}

\date{\today}
\begin{document}
\maketitle
 \pagenumbering{arabic}
 \begin{abstract}
   In industry, shape optimization problems are of utter importance when designing structures such as aircraft, automobiles and turbines.
   For many of these applications, the structure changes over time, with a prescribed or non-prescribed movement.
   Therefore, it is important to capture these features in simulations when optimizing the design of the structure.
   Using gradient based algorithms, deriving the shape derivative manually can become very complex and error prone, especially in the case of time-dependent non-linear partial differential equations.
   To ease this burden, we present a high-level algorithmic differentiation tool that automatically computes first and second order shape derivatives for partial differential equations posed in the finite element frameworks FEniCS and Firedrake.
   The first order shape derivatives are computed using the adjoint method, while the second order shape derivatives are computed using a combination of the tangent linear method and the adjoint method.
   The adjoint and tangent linear equations are symbolically derived for any sequence of variational forms.
   As a consequence our methodology works for a wide range of PDE problems and is discretely consistent.
   We illustrate the generality of our framework by presenting several examples, spanning the range of linear, non-linear and time-dependent PDEs for both stationary and transient domains.
 \end{abstract}
\section{Introduction}
Shape optimization problems constrained by partial differential equations (PDEs) occur in various scientific and industrial applications, for instance when designing aerodynamic aircraft~\cite{reuther1996aerodynamic} and automobiles~\cite{muyl2004hybrid,harbrecht2016}, or acoustic horns~\cite{bangtsson2003shape}.

Shape optimization problems for systems modeled by non-cylindrical evolution PDEs are encountered for example, in fluid-structure and free boundary problems~\cite{moubachir2006moving}.
This category of optimization problems was introduced in the 1970's~\cite{Cea1973} and studied with perturbation theory in the setting of fluid mechanics~\cite{pironneau1974optimum}.
The theoretical foundation for this class of problems where laid by the perturbation of identity method~\cite{murat1975etude}, and the speed method~\cite{zolesio1979identification}.
An recent overview of the development and theory in shape analysis for moving domains can be found in~\cite{moubachir2006moving}.

Mathematically, these problems can be written in the form:
\begin{linenomath}
  \begin{subequations} \label{eq:shape_opt}
\begin{align}
  \min_{\Omega} J(\Omega, u \label{eq:shape_opt:func}),
\end{align}
where $J$ is the objective functional and $u$ is the solution of a PDE over the domain $\Omega$:
\begin{align}
 E(u)= 0 \quad \text{in } \Omega.
\end{align}
\end{subequations}
\end{linenomath}
Here, $E$ denotes the PDE operator. We allow $\Omega$ to be steady (i.e. a static domain) or be time-dependent (i.e. a morphing domain).

Problem \eqref{eq:shape_opt} is typically solved numerically by employing gradient based optimization algorithms,
which require the shape derivative of the goal functional \eqref{eq:shape_opt:func} with respect to the domain $\Omega$.
To obtain an overall fast optimization solver, it is critical that the computation of the shape derivative is efficient.
A finite difference approximation of the shape gradient scales linearly with the number of shape parameters (typically the mesh coordinates), making it in-feasible for many practical problems. The adjoint method is a much more efficient alternative: it yields the first order derivative (the shape gradient) at the equivalent cost of solving one linearized PDE, the adjoint PDE, independent of the number of shape parameters.

Manually deriving and implementing the adjoint and shape derivatives is a laborious and difficult task, especially for time-dependent or non-linear PDEs~\cite{naumann2012art}. Algorithmic differentiation (AD) aims to automate this process by building a computational graph of all elementary mathematical operations in the PDE model. Since the derivative/adjoint of each elementary operation is known, the AD tool can apply the chain rule repeatedly to obtain the full derivative.
This idea has been successfully applied in the context of shape optimization, for instance to the finite volume solvers MIT GCM~\cite{HEIMBACH2005MITGCM}, OpenFOAM~\cite{TOWARA2013429}, SU2~\cite{economon2015su2,zhou2015discrete,sagebaum2017high} and TAU~\cite{gauger2008tau}.
In these works, the AD tool was applied directly to the Fortran or C++ implementation.
One downside of such a `low-level' approach is that the mathematical structure of the forward problem gets intertwined with implementation details, such as parallelization and linear algebra routines, resulting in high memory requirements and a slow-down of 2-10x compared to the theoretical optimal performance~\cite{TOWARA2013429}.

To avoid this intertwining, \cite{funke2013framework,farrell2013automated} introduced a high-level AD framework for models that solve PDEs with the finite element method. The idea is to treat each variational problem in the model as a single operation~\cite{christianson1997giving} (instead of a sequence of elementary linear algebra operations, such as sums, products, etc. as low level AD would do).
The AD tool dolfin-adjoint~\cite{farrell2013automated,MituschEtAl2019} implements this idea within the FEniCS~\cite{LoggMardalEtAl2011} and Firedrake~\cite{rathgeber2017firedrake} frameworks.
The derivative and adjoint of a variational problem (needed by the AD tool) are available through the Unified Form Language (UFL)~\cite{alnaes2014unified} which expresses variational problems symbolically and allows for automatic symbolic manipulation.
The high-level AD approach has a number of advantages compared to low-level AD, for instances near optimal performance and natural parallel support~\cite{sagebaum2018algorithmic}. However, differentiating with respect to the mesh has not been possible in dolfin-adjoint.

The key contribution of this paper is to extend the high-level AD framework in dolfin-adjoint to support shape derivatives for PDE models written in FEniCS/Firedrake. To achieve this, dolfin-adjoint tracks changes in the mesh as part of its computational graph, by overloading the Mesh-class in FEniCS/Firedrake and the corresponding assemble and solve routines.
The shape derivatives of individual variational forms are obtained using pull back to the reference element and G\^{a}teux derivatives, which was recently added to UFL~\cite{ham2019auto}.
We demonstrate that our approach inherits the advantages of dolfin-adjoint with minimal changes to the forward problem, supports first and second order shape derivatives and supports both static and time-dependent domains shapes.

The paper is organised as follows: First, in \cref{sec:tubes}, we give a brief introduction to shape analysis for continuous problems with time-dependent domains.
Then, in \cref{sec:discretetubes,sec:discshape}, we present the same analysis from a discrete shape analysis setting, using finite differences for temporal discretization, and finite elements for spatial discretization.
With this analysis at hand, in \cref{sec:adshop}, we present how to use high-level algorithmic differentiation to differentiate the discrete shape optimization problem.
In \cref{sec:implementation}, we explain which FEniCS/Firedrake operators that had to be overloaded to enable shape sensitivities.
In \cref{sec:FEniCS}, we verify the implementation through a documented example, where we compute the first and second order shape sensitivities of a time-dependent PDE over a morphing domain and verify them with a Taylor convergence study.
In \cref{sec:stokes}, we highlight the generality of the implementation, using different optimization methods to solve an optimization problem with an analytic solution~\cite{pironneau1974optimum}.
Then, in \cref{sec:dfg-3} we compute and verify first and second order shape sensitivities for a time-dependent, non-linear partial differential equation.
Finally, we summarize the core results and findings in \cref{sec:conclusion}.

\section{High-Level AD for shape optimization}\label{sec:ShapeAD}
The goal of this section is to derive a high-level algorithmic differentiation method for the computation of shape derivatives for functionals defined on time-dependent domains.

First, in \cref{sec:tubes} we give a brief summary of the field of continuous shape analysis for time-dependent domains.
Secondly, in \cref{sec:discretetubes}, we consider discrete shape derivatives, where the goal functional on the time-dependent domain is discretized with a finite difference temporal discretization.
Then, in \cref{sec:discshape}, we present the discrete time-dependent shape optimization problem.
Following, in \cref{sec:adshop} we describe how to obtain discretely consistent shape sensitivities using algorithmic differentiation.
Finally, in \cref{sec:gen}, we mention some of the generalizations that is not covered by the previous sections.

\subsection{Continuous shape-analysis on time-dependent domains }\label{sec:tubes}
This subsection gives a brief overview of the main results for shape differentiation over time-dependent domains.
A thorough overview can be found in~\cite{moubachir2006moving}.

\begin{figure}[!ht]
  \includegraphics[width=\linewidth]{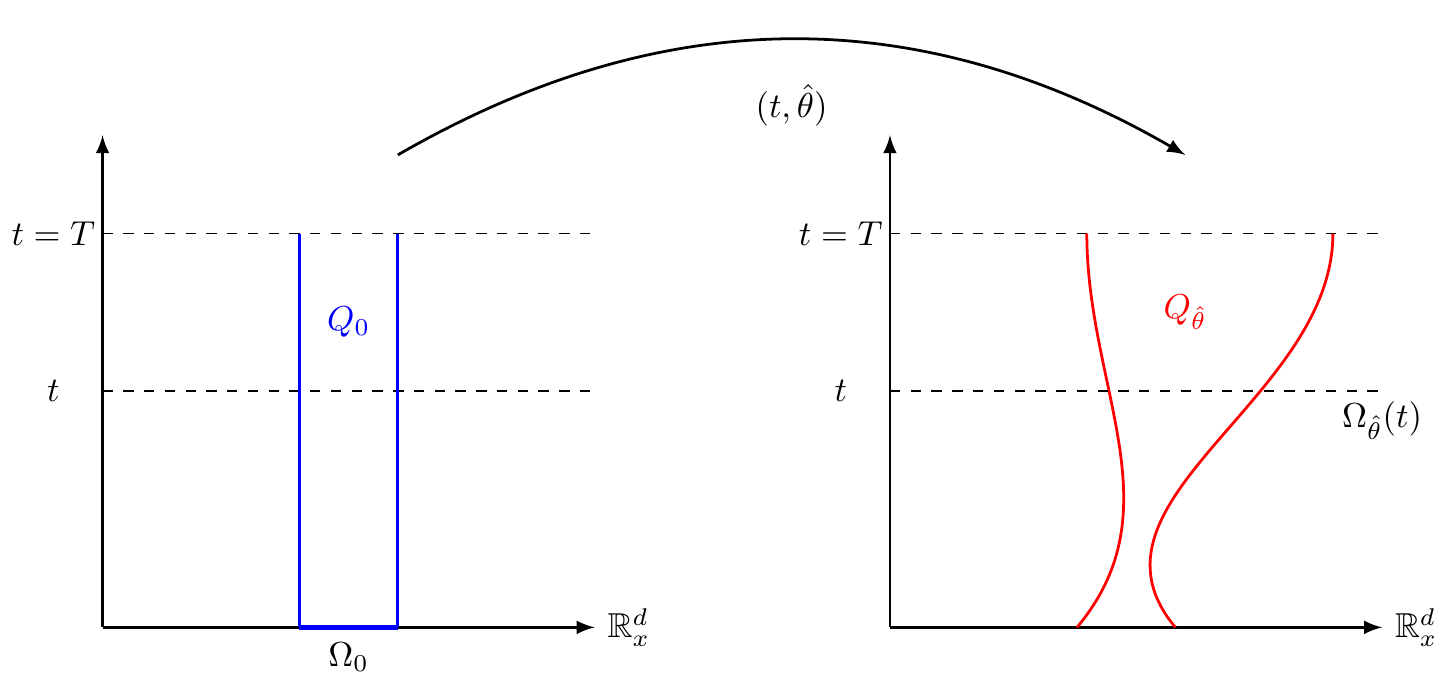}
	\caption{Graphical illustration of a tube $Q_0$ perturbed with the vector field $\hat\theta$.}\label{fig:tube}
\end{figure}

For an initial space domain $\Omega_0\subset\mathbb{R}^d$ we define the space-time tuple $Q_0=[0,T]\times\Omega_0$.
Then, we define a smooth perturbation map $\hat\theta:Q_0\rightarrow\mathbb{R}^d$ and the image
$$\Omega_{\hat\theta}(t):=\hat\theta(t,\Omega_0), \quad \text{for } t\in[0,T].$$
The perturbed, noncylindrical evolution domain is called a tube and defined as
\begin{linenomath}
  \begin{align}
    Q_{\hat\theta} := \cup_{0<t<T} ({t} \times \Omega_{\hat\theta}(t) ).
  \end{align}
\end{linenomath}
The map of the evolution domain is exemplified in \cref{fig:tube}.
We call $\Omega_0$ the base of the tube.

In the setting of shape optimization, we would like to solve the following problem
\begin{linenomath}
\begin{align}
  \min_{Q\in\mathcal{A}} J(Q),
	\label{eq:shape_opt_problem_general}
\end{align}
\end{linenomath}
where $\mathcal{A}$ is the collection of all admissible shape evolution sets.
The cost functional $J$ is typically expressed in terms of integrals over the noncylindrical evolution domain and/or its lateral boundary.
To be able to employ standard differential calculus, we reformulate the functional in terms of $\theta$, that is $j(\hat\theta):=J(Q_{\hat\theta})$.
This formulation is also relevant in cases where the functional is not only dependent on the tube, but the vector field that builds the tube.
This occurs for instance in fluid structure interaction problems and convection diffusion problems, as illustrated in \cref{sec:FEniCS}.
With this formulation, Problem \eqref{eq:shape_opt_problem_general} becomes
\begin{linenomath}
\begin{align}
	\min_{\hat \theta} j(\hat\theta) \quad \text{subject to } Q_{\hat\theta} \in\mathcal{A}. \label{eq:shape_opt_problem_rw}
\end{align}
\end{linenomath}
For the shape analysis of \eqref{eq:shape_opt_problem_rw}, we consider a prototypical example where $j$ is a volume integral:
\begin{linenomath}
  \begin{align}
    j(\hat\theta) &:= \int_{Q_{\hat\theta}} u(\hat\theta)\md Q_{\hat\theta}
    =\int_0^T \int_{\Omega_{\hat\theta}(t)} u(\hat\theta)(t, x) \md x \md t,
    \label{eq:Jproto}
  \end{align}
\end{linenomath}
where $u(\hat\theta)\in H(Q_{\hat\theta})$ is a function in a suitable Hilbert space $H$.

We define the non-cylindrical material derivative of $u$ at $\hat\theta$, $Q_{\hat\theta}\in\mathcal{A}$, in direction $\delta\hat\theta:=\delta\theta \circ\hat\theta$ with $\delta \theta : \mathbb R^d \to \mathbb R^d$ as
\begin{linenomath}
\begin{align}
  \dot u(\hat\theta)\cdot\delta\theta := \lim_{\rho\to0^+} \frac{u((I+\rho\delta\theta)\circ \hat\theta)-u(\hat\theta)}{\rho},
\end{align}
\end{linenomath}
if the limit exists a.e. for $(t,x)\in Q_{\hat\theta}$.
Here $I:\mathbb{R}^d\rightarrow\mathbb{R}^d$ is the identity operator.

The non-cylindrical shape derivative is related to the material derivative by
\begin{linenomath}
\begin{align}
  u'(\hat\theta)\cdot \delta \theta = \dot u(\hat\theta)\cdot \delta \theta - \nabla u (\hat\theta)\cdot\delta \theta.
\end{align}
\end{linenomath}

The method of mappings is used to obtain the shape derivative of \cref{eq:Jproto}. We recall the theorem for the tube derivative of volume functionals, as presented in Chapter 6 of \cite{moubachir2006moving}.
\begin{theorem}[Tube derivative of a volume functional\cite{moubachir2006moving}]\label{thrm:tube}
  For a bounded domain $\Omega_0$ assume that
  $\delta\hat\theta=\delta \theta\circ\hat\theta$ and its inverse is a $\mathcal{C}^1$ differentiable function with respect to all inputs and outputs.
    Then if $u(\hat\theta)$ admits a non-cylindrical material derivative
  $\dot u(\hat\theta)\cdot\delta \theta$, then
    $j(\hat \theta)=\int_{Q_{\hat\theta}}u(\hat\theta)\md Q_{\hat\theta}$ is differentiable at $\hat\theta$ if $\hat\theta$ is in the set of admissible functions, and the derivative of $j$ is given by
  \begin{linenomath}
    \begin{align}
      \begin{split}
    \totder{}{\hat\theta}j(\hat\theta)\cdot\delta\theta &= \int_{Q_{\hat\theta}}\left(\dot u(\hat\theta) + u(\hat\theta)\mathrm{div}\delta\theta \right)\mathrm{d}Q_{\hat\theta}\\
    &= \int_{0}^T\int_{\Omega_{\hat\theta}(t)}\Big(\dot u(\hat\theta)(t, x) + u(\hat\theta)(t, x)\mathrm{div}\delta\theta(t,x)\Big)\md x \md t.
      \end{split}
  \end{align}
  \end{linenomath}
\end{theorem}
A similar result can be derived for functionals involving boundary integrals, see~\cite{moubachir2006moving}.
Theorem \ref{thrm:tube} holds for continuous tubes, that is tubes with a continuous perturbation field $\hat\theta(t,x)$.

To numerically solve partial differential equations and the corresponding shape optimization problem, the tube has to be discretized.
In the optimization community, there are two different pathways to compute sensitivities: first optimize then discretize, and first discretize then optimize.
In this paper, we will consider the discretize-then-optimize strategy, therefore we next consider the temporal discretization of a time-dependent shape optimization problem.

\subsection{Tube derivatives on time-discretized domains}\label{sec:discretetubes}
\begin{figure}[!ht]
  \centering
  \includegraphics[width=0.5\linewidth]{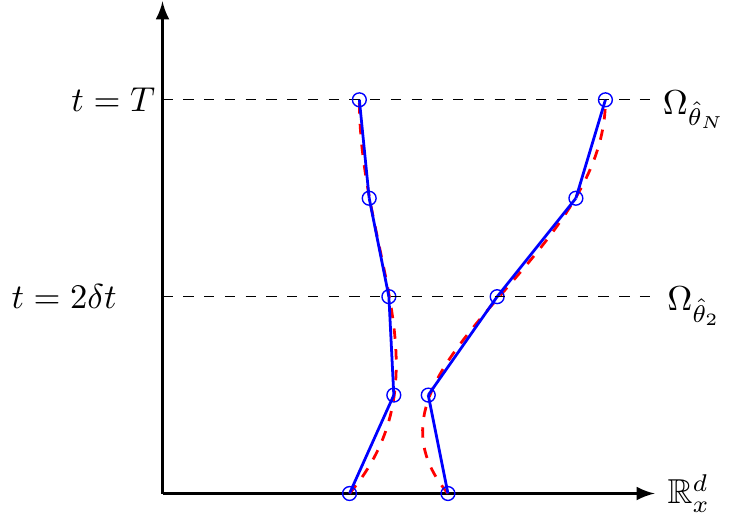}
  \caption{Discretized tube, where each discrete time-step defines an independent perturbation function $\hat\theta_i(x)$ that describes the domain perturbation to the domain at the next time-step.}\label{fig:discretetube}
\end{figure}

For the temporal discretization, we divide the time domain $[0,T]$ into $N$ intervals, separated at $t_i, i=0,\dots,N$.
We also define $N+1$ vector fields $\hat\theta_i(x):\mathbb{R}^d\to\mathbb{R}^d$ that describe the domain perturbations from the $i$th to $(i+1)$th time-step, as visualized in \cref{fig:discretetube}.

We define $\Omega_{\hat\theta_i}:=\hat\theta_i(\Omega_{\hat\theta_{i-1}})=\{x+\hat\theta_i(x) : x\in\Omega_{\hat\theta_{i-1}} \}$ as the $i$th discrete domain used for each time-step, where $\Omega_{\hat\theta_0}=\hat\theta_{0}(\Omega_0)$.
As in \cref{sec:tubes}, we will use the prototypical functional \eqref{eq:Jproto} to illustrate the concepts of discretized tube derivatives.
We note at this point that the algorithmic differentiation (AD) framework presented in \cref{sec:adshop} generalizes to a wide range of functionals, including boundary integrals, products of integrals etc.

We use a generalized finite difference scheme to rewrite \eqref{eq:Jproto} and obtain the time-discretized functional
\begin{linenomath}
\begin{align}
  j(u, \hat\theta) := \sum_{i=0}^Nw_ij_i(u_i, \hat \theta_i)=\sum_{i=0}^Nw_i \int_{\Omega_{\hat\theta_i}}u(\hat\theta_i)(x)\md x,
\end{align}
\end{linenomath}
where $w_i$ is the $i$th finite difference weight and $\hat \theta := (\hat \theta_0, ..., \hat \theta_N)$.
Following the same steps as in \cref{thrm:tube}, one obtains the shape gradient of the time-discretized functional:
\begin{linenomath}
  \begin{align}\label{eq:disc:grad}
  \totder{}{\hat\theta}j(u, \hat\theta)\cdot\delta\theta
  &= \sum_{i=0}^Nw_i \int_{\Omega_{\hat\theta_i}}\Big(\dot u(\hat\theta_i)(x)\cdot\delta\theta_i(x) + u(\hat\theta_i)(x)\mathrm{div}\delta\theta_i(x)\Big)\md x,
\end{align}
 \end{linenomath}
where $\delta\theta_i$ is the perturbation function at the $i$th time-step, and the material derivative $\dot u(\hat\theta_i)$ is defined as
\begin{linenomath}
  \begin{align}
    \dot u(\hat\theta_i):=\lim_{\rho\to 0} \frac{u((I+\rho\delta\theta_i)\circ \hat\theta_i)-u(\hat\theta_i)}{\rho}.
  \end{align}
\end{linenomath}

\subsection{Discrete time-dependent shape optimization problems with PDE constraints}\label{sec:discshape}
Using the notation from \cref{sec:tubes}, we can write a continuous shape optimization problem with PDE constraints as
\begin{linenomath}
  \begin{subequations}\label{eq:StrongOpt}
  \begin{align}
  \min_{u, \hat\theta} j(u,\hat\theta) \\
  \text{ subject to }
  E(u, \hat\theta)&=0 \quad \text{in } Q_{\hat\theta},
  \end{align}
  \end{subequations}
\end{linenomath}
where $j(u, \hat \theta)$ is the (non-discretized) goal functional, $E(u, \hat\theta)$ is a time-dependent PDE operator with solution $u$ over $Q_{\hat\theta}$, where $Q_{\hat\theta}$ is the non-cylindrical evolution domain, as described in \cref{sec:tubes}.

As in \cref{sec:discretetubes}, we discretize the problem in time with finite differences. This yields a sequence $E_0(u_0),\cdots, E_N(u_0,\cdots,u_N)$ of PDE operators for each time-step.

Each PDE operator is discretized in space using the finite element method (FEM) for finding a numerical approximation to the PDE.
Therefore, we find the variational formulation for each PDE operator $E_i$ by multiplying with a test-function $v\in V_i(\Omega_{\hat\theta_i})$, and performing integration by parts if needed. We use $F_i$ to denote the corresponding variational formulation of $E_i$.

With that, the discretized version of problem \eqref{eq:StrongOpt} reads
\begin{linenomath}
  \begin{subequations}\label{eq:DiscOpt}
\begin{align}
  &\min_{\hat\theta_0,\dots,\hat\theta_N}\sum_{i=0}^Nw_ij_i(u_i, \hat \theta_i) \label{eq:disc:opt}
\end{align}
where $\bar{u}_0,\cdots,\bar{u}_N$ are the implicit solution operators solving
  \begin{align}
    \begin{split}\label{eq:disc:varprob2}
      F_0(u_0,\hat\theta_0; v)& = 0\quad\forall v\in V_0(\Omega_{\hat\theta_0}),\\
      & \vdots\\
      F_N(u_0,\hat\theta_0,\cdots, u_N,\hat\theta_N; v)&=0 \quad\forall v\in V_N(\Omega_{\hat\theta_N}).
    \end{split}
  \end{align}
  \end{subequations}
\end{linenomath}
If at one time-step the domain changes in the inwards normal direction, then there exist points $x\in \Omega_{\hat\theta_{i-1}}$ where $x\notin \Omega_{\hat\theta_{i}}$, and thus previous solutions must be mapped to $\Omega_{\hat\theta_{i}}$. Examples of such mappings in the continuous and finite element setting are given in~\cite{berggren}.




\subsection{Algorithmic differentiation for the discrete shape optimization problem}\label{sec:adshop}
This subsection explains how to compute the discretely consistent shape gradient for Problem \eqref{eq:DiscOpt} using algorithmic differentiation (AD).
The fundamental idea of algorithmic differentiation is to break down a complicated numerical computation, like a numerical finite element model, into a sequence of simpler operations with known derivatives.
By systematic application of the chain rule, one obtains the derivative of the composite function using only partial derivatives of these simple operations~\cite{griewank2008}.

There are two different modes of AD, namely the \textit{forward mode} and \textit{reverse mode}.
Forward mode AD computes directional derivatives, while reverse mode AD computes gradients.
Hence, forward mode AD is most often applied when the number of outputs are greater than the number inputs, and the reverse mode is used in the opposite case.

In this paper, we consider first order meshes, thus meshes where each cell is defined by their vertices.
Therefore, the discrete control variable will be a vector-function with degrees of freedom on each of the vertices.
This corresponds to a function in the finite element function space of piecewise continuous, element linear functions.
Typically, such a function has thousands to millions of degrees of freedom.
Therefore the reverse mode is a popular choice for first order derivatives.
For second order derivatives, a combination of forward and reverse mode is the most efficient \cite{naumann2012art}.

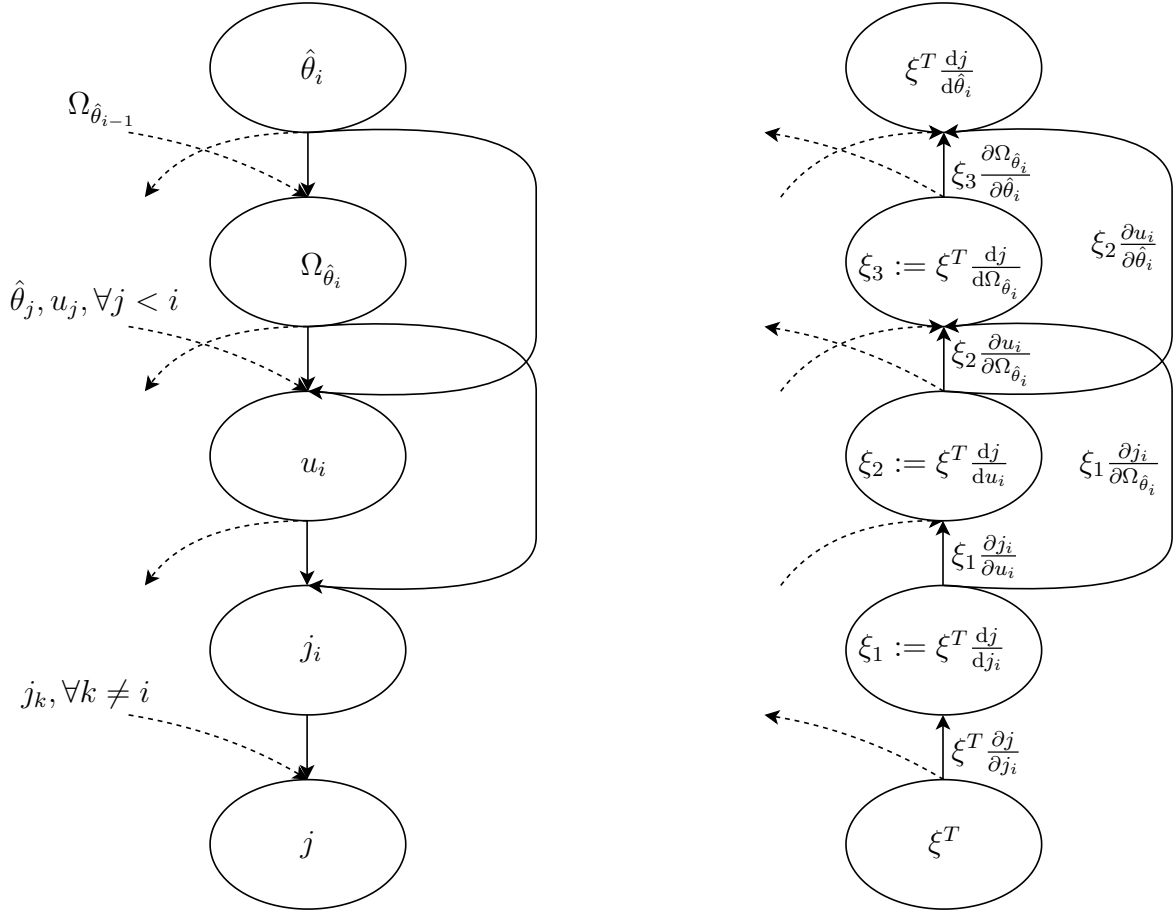
\begin{figure}[!ht]
  \centering
\vspace{1cm}
\begin{overpic}[width=0.9\linewidth,grid=false]{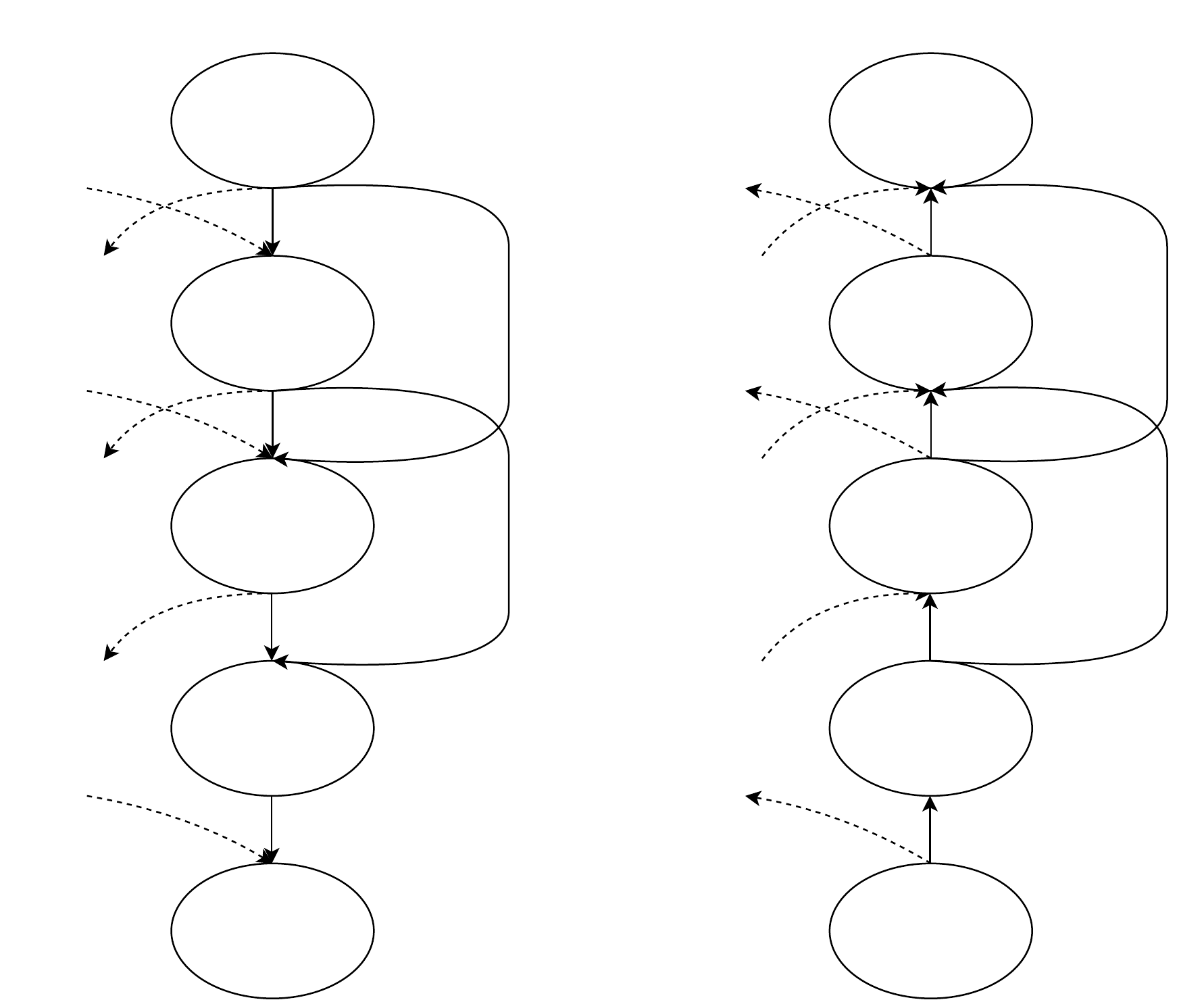}
    \put(2,69){$\Omega_{\hat\theta_{i-1}}$}
    \put(22,72){$\hat\theta_{i}$}
    \put(22,55){$\Omega_{\hat\theta_{i}}$}
    \put(22,38){$u_i$}
    \put(22,22){$j_i$}
    \put(22,5){$j$}
    \put(-3,52){$\hat\theta_j,u_j, \forall j<i$}
    \put(-2,18){$j_k,\forall k\neq i$}
    \put(74,72){$\xi^T \totder{j}{\hat\theta_{i}}$}
    \put(70,55){$\xi_3 := \xi^T\totder{j}{\Omega_{\hat\theta_{i}}}$}
    \put(70,38){$\xi_2 := \xi^T\totder{j}{u_i}$}
    \put(70,22){$\xi_1 := \xi^T\totder{j}{j_i}$}
    \put(76,5){$\xi^T$}
    \put(78,13){$\xi^T\der{j}{j_i}$}
    \put(78,30){$\xi_1\der{j_i}{u_i}$}
    \put(78,47.5){$\xi_2\der{u_i}{\Omega_{\hat \theta_{i}}}$}
    \put(89,38){$\xi_1\der{j_i}{\Omega_{\hat\theta_i}}$}
    \put(90,57){$\xi_2\der{u_i}{\hat\theta_{i}}$}
    \put(78,63){$\xi_3\der{\Omega_{\hat\theta_i}}{\hat\theta_{i}}$}
\end{overpic}
  \caption{The acyclic subgraph for the operations required on the $i$th time-step of the partial differential equation.
  On the left hand side the forward computational subgraph is illustrated,
  with the resulting reverse AD graph on the right hand side.
  The initial $\xi$ scales the derivative and is usually chosen to be $1$.}\label{fig:graph}
\end{figure}

In order to apply AD to the discretized shape optimization problem \eqref{eq:DiscOpt}, we decompose our model into four unique operations:
\begin{enumerate}
\item Domain perturbation at the $i$th time-step:
\begin{equation}
\Omega_{\hat\theta_{i}} = \hat\theta_i(\Omega_{\hat\theta_{i-1}}).\label{eq:op1}
\end{equation}
\item PDE solver at the $i$th time-step, solving the variational problem \eqref{eq:disc:varprob2}:
\begin{equation}
u_i = \bar u_i(u_0,\hat\theta_0,\cdots, u_{i-1},\hat\theta_i, \Omega_{\hat\theta_i}).\label{eq:op2}
\end{equation}
\item Spatial integration of the functional at the $i$th time-step:
\begin{equation}
      j_i = \int_{\Omega_{\hat\theta_i}}u_i\md x.\label{eq:op3}
\end{equation}
\item Temporal integration of the functional:
\begin{equation}
      j = \sum_{i=0}^N w_i j_i.\label{eq:op4}
\end{equation}
\end{enumerate}
With these operations, we can create a computational graph for the functional evaluation of \eqref{eq:disc:opt}.
The left side of \cref{fig:graph} illustrates the subgraph associated with $j_i$.
The edges represent dependencies between variables, where the arrows are pointing in the direction information is flowing.
The dashed lines represent any number of dependencies which enter from or exit the subgraph, i.e.\ edges that connect with a subgraph for $j_k$ with $k \neq i$.
A node is illustrated as an ellipse.
We denote nodes without incoming edges as \textit{root nodes}.

In forward mode AD, the labels in the forward graph can be substituted with directions or perturbations.
For non-root nodes, these directions are computed as the partial derivative in the direction of predecessor nodes.
For example, the direction at the ellipse for $\Omega_{\hat\theta_i}$ will be
\begin{linenomath}
\begin{align*}
    \delta \Omega_{\hat\theta_i} := \der{\Omega_{\hat\theta_i}}{(\Omega_{\hat\theta_{i-1}}, \hat\theta_i)}[\delta \Omega_{\hat\theta_{i-1}}, \delta \hat\theta_i],
\end{align*}
\end{linenomath}
where $\delta \Omega_{\hat\theta_{i-1}}$ and $\delta \hat\theta_i$ are the directions at the $\Omega_{\hat\theta_{i-1}}$ and $\delta \hat\theta_i$ ellipses, respectively.
In the case of root nodes, the user specifies some initial direction.

When performing reverse mode AD, the flow of information is reversed.
Thus, the computational graph is reversed with arrows pointing in the opposite direction and new operations are associated with the edges and nodes.
To start off a reverse mode AD, a weight $\xi$ in the codomain of the forward functional $j$ is chosen.
The weight can be thought of as a vector with the result of the AD computations being the vector-Jacobian product $\xi^T Dj$, where $Dj$ is the Jacobian matrix of $j$.

The right side of \cref{fig:graph} illustrates the reverse mode AD for the subgraph of $j_i$.
Each node in the reverse graph is associated with the corresponding node in the forward graph.
The last node in the forward graph is associated with the first node of the reverse graph etc.
Unlike the forward graph, the edges now represent the propagation of a different variable than the one found inside the ellipses.
An outgoing edge from an ellipse represents the product of the value inside the ellipse and the partial derivative of the associated
variable of the upstream node, with respect to the variable associated with the downstream node.
For example, the edge in the reverse AD graph associated with $j$ and $j_i$ represents the product $\xi^T \der{j}{j_i}$.
At the points where multiple arrows meet the values of the edges are summed producing the result inside the ellipse.
Thus, the values inside the ellipse is the gradient or total derivative of $\xi^T j$ with respect to the variable associated with the ellipse in the forward graph.
For brevity, the values along the dashed lines are omitted.

For each forward operation, the AD tool needs to have access to the partial
derivatives with respect to its dependencies. In forward AD mode, the partial
derivative is multiplied from the right with a direction $\delta$, while in
reverse AD mode it is left multiplied with a weight $\xi$.  We will now go
through the operations in the order they are encountered in the reverse mode.

\subsubsection{The summation operator $j$}
\label{sec:ad:operations:sumj}
Considering the operation $j=\sum_{i=0}^Nw_ij_i$, the partial derivative with respect to $j_i$ is $\der{j}{j_i}=w_i$.
For a direction $\delta$ or weight $\xi$, the right and left multiplications are $w_i \delta$ and $\xi w_i$, respectively.

\subsubsection{The integral operator $j_i$}\label{sec:opt3}
Next we consider the operation $j_i(u_i, \hat\theta_i) = \int_{\Omega_{\hat\theta_i}}u_i\md x$.
The partial derivative with respect to $u_i$ 
can directly be obtained using standard differentiation rules.
The partial derivative $\der{j_i}{\Omega_{\hat\theta_i}}$ is slightly more complicated.
To obtain a discretely consistent partial shape derivative for a functional containing finite element functions, we do a brief recollection of the core results of~\cite{ham2019auto}.

Let $\{K_l\}_{l\in\mathcal{L}}$ be a partition of $\Omega_{\hat\theta_i}$ such that the elements $K_l$ are non-overlapping, and $\cup_l\overline{K}_l=\overline{\Omega(\hat\theta_i)}$.
We denote the mapping from the reference cell $\hat{K}$ to $K_l$ as $\phi_l(\hat K)$,
for each $l\in\mathcal{L}$.
Consider the perturbation function $\tau^i_\rho(x) = x+ \rho \delta\theta_i(x), \rho\in[0,\alpha], x\in\Omega(\hat\theta_i)$. Thus, the perturbed domain can be written as the partition $\{\tau^i_\rho(K_l)\}_{l\in\mathcal{L}}=\{\tau^i_\rho\circ\phi_l(\hat K)\}_{l\in\mathcal{L}}$.
Using the finite element discretization and change of variables, we rewrite the integral operation \eqref{eq:op3} as an integral over the reference element
\begin{linenomath}
\begin{align}
  j_i(u_i, \rho\delta\theta_i)&=\sum_{l\in\mathcal{L}}\Int{\tau^i_\rho(K_l)}{}u_i(\tau^i_\rho(x))\md x
  =\sum_{l\in\mathcal{L}}\Int{\hat{K}}{}(u_i\circ\tau^i_\rho\circ\phi_l)
  \vert\mathrm{det}(D(\tau^i_\rho\circ \phi_l))\vert\md x.
\end{align}
\end{linenomath}
As shown in~\cite{ham2019auto} the shape derivative can be written using the G\^ateaux derivative of the map $T\mapsto ((u_i)_{T\circ \inv{\phi_l}}\circ T \mathrm{det}(DT))$ at $T=\phi_l$ in direction $\delta\theta_i\circ \phi_l$:
\begin{linenomath}
\begin{align}\label{eq:functional_shape_derivative}
  \der{j_i}{\Omega_{\hat\theta_i}}[\delta\theta_i]&= \sum_{l\in\mathcal{L}}\Int{\hat{K}}{}
  \langle d_T\left[((u_i)_{ T\circ\phi_l^{-1}}\circ T)\vert \mathrm{det}(DT)\vert\right],
    \delta\theta_i\circ \phi_l\rangle\vert_{T=\phi_l}\md x,
\end{align}
\end{linenomath}
which is the directional derivative with direction $\delta \theta_i$ required for forward mode AD.

When performing reverse mode AD with a weight $\xi$, the full Jacobian of \eqref{eq:functional_shape_derivative} assembled and multipled by $\xi$.
The derivation of $\der{j_i}{\Omega_{\hat\theta_i}}[\cdot]$ is automatically computed in FEniCS using~\cite{ham2019auto}.

\subsubsection{The implicit PDE operation $\bar u_i$}
For the implicit function $\bar u_i$, the output $u_i\in W_i(\Omega_{\hat\theta_i})$ is the result of the relation
\begin{linenomath}
  \begin{equation}
    F_i(u_0,\hat\theta_0,\cdots, u_i,\hat\theta_i,\Omega_{\hat\theta_i}; v)=0 \quad\forall v\in V_i(\Omega_{\hat\theta_i}).
  \end{equation}
\end{linenomath}
Let us consider a placeholder variable $m\in M$, where $M$ is the appropriate vector space,
which could be any of the dependencies of $\bar u_i$.
For the PDE solution of $u_i$, we require two operations $\der{\bar{u}_i}{m}\delta_m$ and $\xi\der{\bar{u}_i}{m}$
where $\delta_m\in M$ and $\xi\in W_i(\Omega_{\hat\theta_i})$ are the results of previous computations of the forward and reverse mode, respectively.

%
For forward mode the directional derivative can be computed using the tangent linear model of the PDE
\begin{align}
\der{\bar u_i}{m}\delta_m = -\der{F_i}{u_i}^{-1}\der{F_i}{m}\delta_m.
\end{align}

In the reverse mode, the derivative is computed in two steps.
First the adjoint equation is solved
\begin{align}
  \der{F_i}{u_i}^*\lambda=\xi,
\end{align}
where the $*$ denotes the Hermitian adjoint.
Second, the partial derivative of $\bar u_i$ with respect to any variable $m$ can be computed as
\begin{align}
  \xi\der{\bar u_i}{m} = -\lambda^*\der{F_i}{m}.
\end{align}


Thus, if $m$ is equal to $\Omega_{\hat\theta_i}$, the derivative is computed on the reference element, as described for $j_i$ in the previous section.

\subsubsection{The domain perturbation operator $\hat\theta_i$}

The mesh perturbation operator $\hat\theta_i(x) = x + \hat\theta_i(x)$ is linear in both $x$ and $\hat\theta_i$, and its derivatives are the identity operations.
Thus, for forward mode AD with direction $\delta$ or reverse mode AD with weight $\xi$, the result is $\delta$ or $\xi$, respectively.

\subsection{Generalizations}\label{sec:gen}
In the previous sections, we considered a prototypical example for the functional $j$, where there were no explicit dependencies of $\hat\theta_i$ in the integrand, and the integrand was not a function of spatial derivatives.
However, as shown in the next section, this is not a limitation of the algorithmic differentiation framework.
Additionally, the previous sections did not explicitly handle boundary conditions.
These can be handled either strongly or weakly in the proposed framework.
\section{Implementation}\label{sec:implementation}
To solve Problem \eqref{eq:DiscOpt} numerically, we use the FEniCS project~\cite{LoggMardalEtAl2011} and dolfin-adjoint\cite{MituschThesis2018}.
The FEniCS project is a framework for solving PDEs using the finite element method. It uses the Unified Form Language~\cite{alnaes2014unified} to represent variational forms in close to mathematical syntax. UFL has support for symbolic differentiation of forms, and recently, shape derivatives, see  Ham et al.~\cite{ham2019auto}. The user-interface of the FEniCS-project is called dolfin~\cite{logg2010dolfin}, and has both a Python and C++ user-interface.
dolfin-adjoint is a high-level algorithmic differentiation software, that uses operator overloading to augment dolfin with derivative operations.
dolfin-adjoint implements both tangent linear (forward) and adjoint (reverse) mode algorithmic differentiation.
Second-order derivatives are implemented using \textit{forward-over-reverse mode}, where tangent linear mode is applied to the adjoint model.

For this paper, we have extended dolfin-adjoint to compute shape derivatives of FEniCS models. The following subsections will go through these extensions.

\subsection{The domain perturbation operator $\hat\theta_i$}
Since the domain perturbation operation has the computational domain, represented by the \pythoninline{dolfin.Mesh}, and a \pythoninline{dolfin.Function} as input, these two classes is overloaded such that they can be added to the computational graph. The operator in dolfin-adjoint which represents $\hat\theta_i$ is the \pythoninline{ALE.move} function. Therefore, we have added the operations required to evaluate first and second order derivatives, as required by the different AD modes.

\subsection{The implicit PDE operation $\bar u_i$}
The simplest way of solving a PDE in FEniCS, is to write the variational formulation in UFL, then call \pythoninline{solve(Fi==0, ui, bcs=bc)}, where \pythoninline{Fi} is the $i$th variational formulation, \pythoninline{ui} the function to solution is written to, and \pythoninline{bcs} a list of the corresponding Dirichlet boundary conditions.
We extended the overloaded solve operator in pyadjoint to differentiate with respect to \pythoninline{dolfin.Mesh}, as explained in \cref{sec:opt3}. 

\subsection{The integral operator $j_i$}
Integration of variational formulations and integrals written in UFL is performed by calling the \pythoninline{assemble}-function.
This function can return a scalar, vector or matrix, depending on the form $j_i$.
This operator has been extended with shape derivatives, as explained in \cref{sec:opt3}.
In general, the \pythoninline{assemble} function can be used in combination with the implicit PDE operation $\bar u_i$, in for instance \pythoninline{KrylovSolver} and \pythoninline{PETScKrylovSolver}.

\subsection{The summation operator}
This operation has been overloaded in pyadjoint, and no additions was required for shape derivatives.

\subsection{Firedrake}
Since Firedrake uses the same high-level user interface to solve PDEs, the \pythoninline{solve} and \pythoninline{assemble} implementation only has minor differences.
However, Firedrake has a unique handling of meshes.
Therefore, the overloading of the mesh class differs from the one used in dolfin.
The mesh perturbation command \pythoninline{ALE.move(mesh, perturbation)}, is replaced by\\
\pythoninline{mesh.coordinates.assign(mesh.coordinates + perturbation)}.
\section{Documented demonstration, application and verification of tube derivatives in FEniCS.}\label{sec:FEniCS}

In this section, we will illustrate how dolfin-adjoint can be used to solve problems with time-dependent domains, highlighting key implementation aspects along the way.

\begin{figure}[!ht]
  \centering
  \begin{overpic}[width=0.35\linewidth]{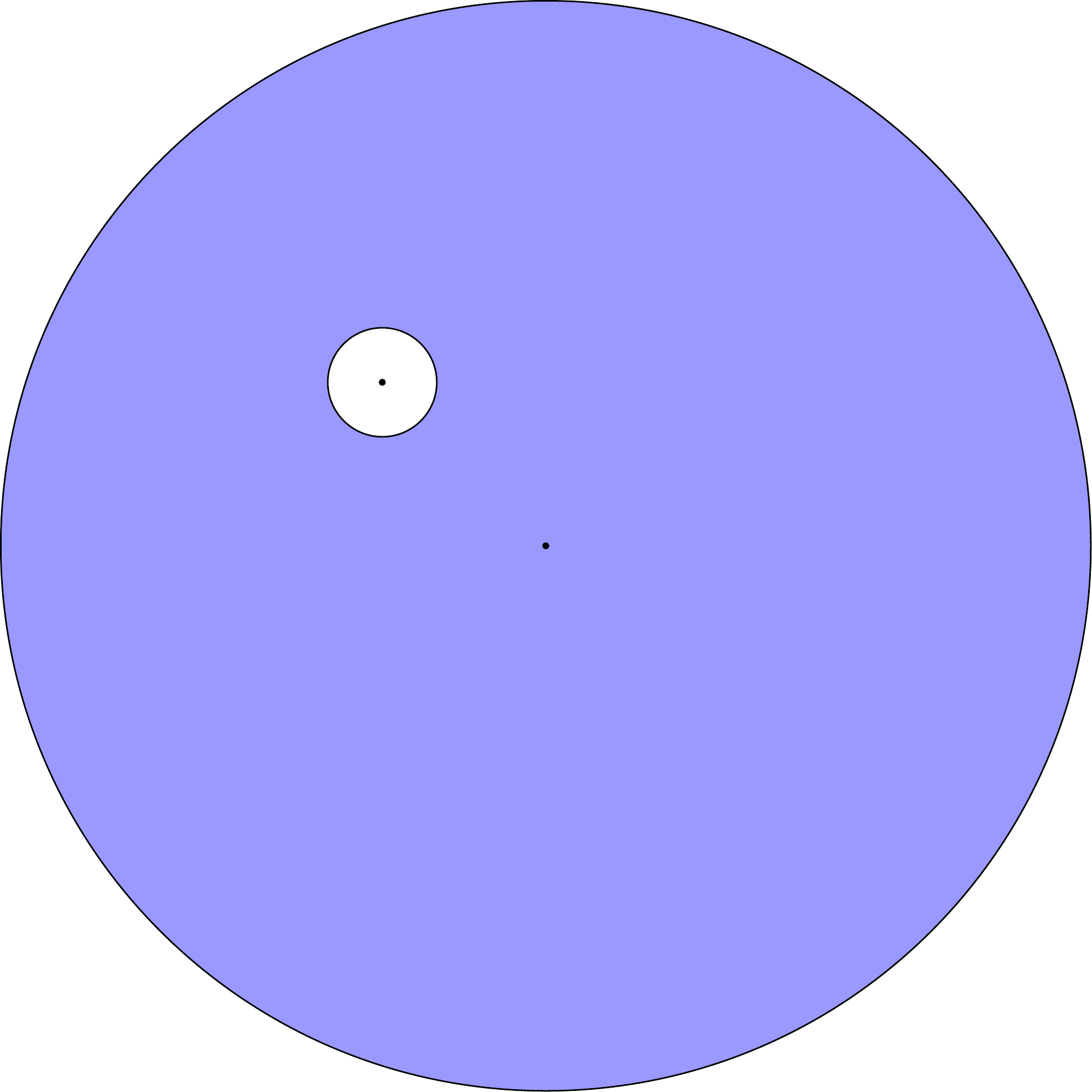}
    \put(40,60){$\partial\Omega_{\hat\theta}^D(0)$}
    \put(80,90){$\partial\Omega_{\hat\theta}^N(0)$}
    \put(50,45){$(0,0)$}
  \end{overpic}
  \caption{Initial domain (Stem of the tube) for \cref{eq:tube:F} with the corresponding boundaries.}\label{fig:init_domain}
\end{figure}

Consider the following problem:
Compute $\totder{J}{\hat\theta}[\delta\theta]$ where
\begin{linenomath}
\begin{align}\label{eq:tube:J}
  J = \int_{Q_{\hat\theta}} \nabla u :\nabla u \md Q_{\hat\theta},
\end{align}
\end{linenomath}
and $u$ is the solution of the advection-diffusion equation
\begin{linenomath}
  \begin{subequations}
\begin{align}\label{eq:tube:F}
  \der{u}{t}-k\Delta u - \nabla \cdot (u \der{\hat\theta}{t}) &= 0 && \text{on } Q_{\hat\theta},\\
  k\der{u}{n}&=0 && \text{on } \partial\Omega_{\hat\theta}^N(t), t\in(0,T),\\
  u&=1 && \text{on } \partial\Omega_{\hat\theta}^D(t), t\in(0,T),\\
  u(x,y,0) &=0 && \text{on } \Omega_{\hat\theta}(0).
\end{align}
  \end{subequations}
\end{linenomath}
The advection velocity $\der{\hat\theta}{t}$ is the time derivative of the domain deformation $\hat\theta$.
In this example, the stem of the tube will be defined by a circular domain with a circular hole, as depicted in \cref{fig:init_domain}.
We choose the initial perturbation velocity field
$\der{\hat\theta}{t}(t,x)=rot(x(t)) = rot(x_1(t),x_2(t))=(2\pi\omega x_2(t), - 2\pi\omega x_1(t))^T$.
The physical interpretation of this setup is that the hole is rotating around the center of the circular domain.

In FEniCS, we start by importing dolfin and dolfin-adjoint, which is overloading core operations of dolfin.
\begin{python}[numbers=left]
  from dolfin import *
  from dolfin_adjoint import *
\end{python}
The next step is to load the discrete representation of the domain, and the facet markers corresponding to markers on the two boundaries of $\Omega_0$.
\begin{python}[numbers=left, firstnumber=3]
  mesh = Mesh()
  with XDMFFile("mesh/mesh.xdmf") as xdmf:
      xdmf.read(mesh)
  mvc = MeshValueCollection("size_t", mesh, 1)
  with XDMFFile("mesh/mf.xdmf") as infile:
      infile.read(mvc, "name_to_read")
  bdy_markers = cpp.mesh.MeshFunctionSizet(mesh, mvc)
\end{python}
Here, dolfin-adjoint overloads the \pythoninline{dolfin.Mesh}-class, as it is the integration domain that is input to the discretized variational formulation.

Next, we define relevant physical quantities and time discretization variables.
\begin{python}[numbers=left, firstnumber=10]
  k = Constant(0.01)     # Diffusion coefficient
  omega = Constant(0.25) # Rotation velocity
  T = 4                  # Final time
  dt = Constant(1e-2)    # Time-step
  N = int(T/float(dt))   # Number of time-steps
\end{python}
To describe the initial mesh movement discretely, we discretize the equation for
$\der{\hat\theta}{t}$ with a Crank-Nicholson scheme in time, yielding:
Find $\hat\theta_n:=x(t_n)-x(t_{n-1})$ such that for all test-functions $z\in V$
\begin{linenomath}
\begin{align}
  (\hat\theta_n, z)_\Omega = \frac{1}{2}\Delta t(rot(x(t_{n-1})+\hat\theta_n) + rot(x(t_{n-1})), z)_\Omega.
\end{align}
\end{linenomath}
In FEniCS, the deformation field is defined as a CG-1 field, where the degrees of freedom are on each of the vertices of the element. The variational form is written in the Unified Form Language~\cite{alnaes2014unified}, yielding
\begin{python}[numbers=left, firstnumber=15]
  V = VectorFunctionSpace(mesh, "CG", 1)
  z = TestFunction(V)
  X = SpatialCoordinate(mesh)
  rot = lambda y: 2*pi*omega*as_vector((y[1], -y[0]))
  F_s = lambda thn: inner(thn, z)*dx\
      - dt*0.5*inner(rot(X+thn)+rot(X), z)*dx
\end{python}
The next step is to discretize \cref{eq:tube:F} with a Crank-Nicholson discretization scheme, yielding: Find $u_n\in W_h^1$ such that for all $v\in W_h^0$
\begin{linenomath}
\begin{align}
  \frac{1}{\Delta t}(u_n-u_{n-1}, v)_\Omega + k (\nabla u_{n-1/2},\nabla v)_{\Omega}
  +\frac{1}{\Delta t}(u_{n-1/2}\hat\theta_{n-1/2}, v)_\Omega&=0,
\end{align}
\end{linenomath}
where $u_{n-1/2}=\frac{1}{2}(u_n+u_{n-1})$, $\hat\theta_{n-1/2}=\frac{1}{2}(\hat\theta_n+\hat\theta_{n-1})$.
We start by creating the variational form symbolically, as it will be re-used for every time-step. We let \pythoninline{F_u} be a function of the mesh velocity \pythoninline{V}.
\begin{python}[numbers=left, firstnumber=21]
  W = FunctionSpace(mesh, "CG", 1)
  u0 = Function(W, name="u^{n-1}")
  v,  w = TestFunction(W), TrialFunction(W)
  F_u = lambda V: (1.0/dt*(w-u0)*v*dx
                    + k*inner(grad(v), 0.5*(grad(w)+grad(u0)))*dx
                    + inner(0.5*(w+u0)*V, grad(v))*dx
\end{python}
Additionally, we create the corresponding Dirichlet condition for the boundary $\partial\Omega_{\hat\theta_i}^D$ which are specified through the facet function \pythoninline{bdy_markers}.
\begin{python}[numbers=left, firstnumber=27]
  bc = DirichletBC(W, Constant(1.0), bdy_markers, 2)
\end{python}
The list of perturbation functions for each time-step is then created with the following command
\begin{python}[numbers=left,firstnumber=28]
  thetas = [Function(V) for i in range(N+1)]
\end{python}
If the initial domain should be controlled, one perturbs the domain
\begin{python}[numbers=left,firstnumber=29]
  ALE.move(mesh, thetas[0])
\end{python}
The function \pythoninline{ALE.move} is an implicit function, perturbing the mesh-coordinates with the CG-1 field \pythoninline{thetas[0]}.
The forward problem in then solved and the functional computed with the for-loop shown in \cref{code:approach:1}.
\begin{python}[caption={The forward simulation of \cref{eq:tube:F}, where the domain movement is input to the control parameter $\hat\theta_i$.},captionpos=b, label={code:approach:1},numbers=left,firstnumber=30]
J = 0
for i in range(N):
    # Solve for mesh displacement.
    with stop_annotating():
        solve(F_s(thetas[i+1])==0, thetas[i+1])

    # Move mesh
    ALE.move(mesh, thetas[i+1])

    # Solve for state
    a, L = system(F_u(0.5/dt*(thetas[i+1]+thetas[i])))
    solve(a==L, u1, bc)
    u0.assign(u1)

    # Compute functional
    J += assemble(dt*inner(grad(u1), grad(u1))*dx)
\end{python}
For each iteration in the for-loop in \cref{code:approach:1}, we obtain a computational sub-graph similar to \cref{fig:graph}.
The first addition to the computational graph is the \pythoninline{ALE.move}-command in line 37.
Then, the \pythoninline{solve} command in line 41 is added to the graph,
and finally the \pythoninline{assemble}-function in line 40 is added to the computational graph.

\subsection{Fixed rotational motion}
The mesh-movement PDE $F_s$ that is solved in \cref{code:approach:1} is not represented in the computational graph, due to the operation \pythoninline{with stop_annotating()}. This means that we only consider rotation as an initial movement for the domain, but that the shape derivative will not restrict changes in the domain to be rotational.

To obtain as system respecting the rotational motion, one can replace line 31-40 in \cref{code:approach:1} with \cref{code:approach:2}, where we have decomposed the perturbation field into a static component, the rotation $S$, and the varying component $\hat\theta_i$.

\begin{python}[caption={The forward simulation of \cref{eq:tube:F}, where the domain movement is decomposed into a fixed component (the rotation) and the control variable component.},captionpos=b, label={code:approach:2}]
# Total deformation field per time-step
S_tot = [Function(V) for i in range(N+1)]
S_tot[0].assign(thetas[0])
for i in range(N):
    # Solve for mesh displacement.
    solve(F_s(S)==0, S)
    S_tot[i+1].assign(S + thetas[i+1])

    # Move mesh
    ALE.move(mesh, S_tot[i+1])

    # Solve for state
    a, L = system(F_u(0.5/dt*(S_tot[i]+S_tot[i+1])))
\end{python}

Due to this change, we obtain additional blocks in the computational graph when solving the rotational system \pythoninline{solve(F_s(S)==0,S)}, and when we assign the two movement vectors to a total vector, \pythoninline{S_tot[i+1].assign(S+thetas[i+1])}.

\subsection{Verification}

With the full code for computing the forward problem, we define the reduced functional $\hat J$, which is a function of the perturbations $\hat\theta_i$.
\begin{python}
 ctrls = [Control(c) for c in thetas]
 Jhat = ReducedFunctional(J, ctrls)
 dJdctrl = Jhat.derivative()
\end{python}
The \pythoninline{Jhat.derivative} call on the last line applies the shape AD framework, solving the corresponding adjoint equation and computing the shape derivatives.

\begin{figure}[!ht]
  \centering
  \begin{subfigure}{\linewidth}
    \includegraphics[width=0.24\linewidth]{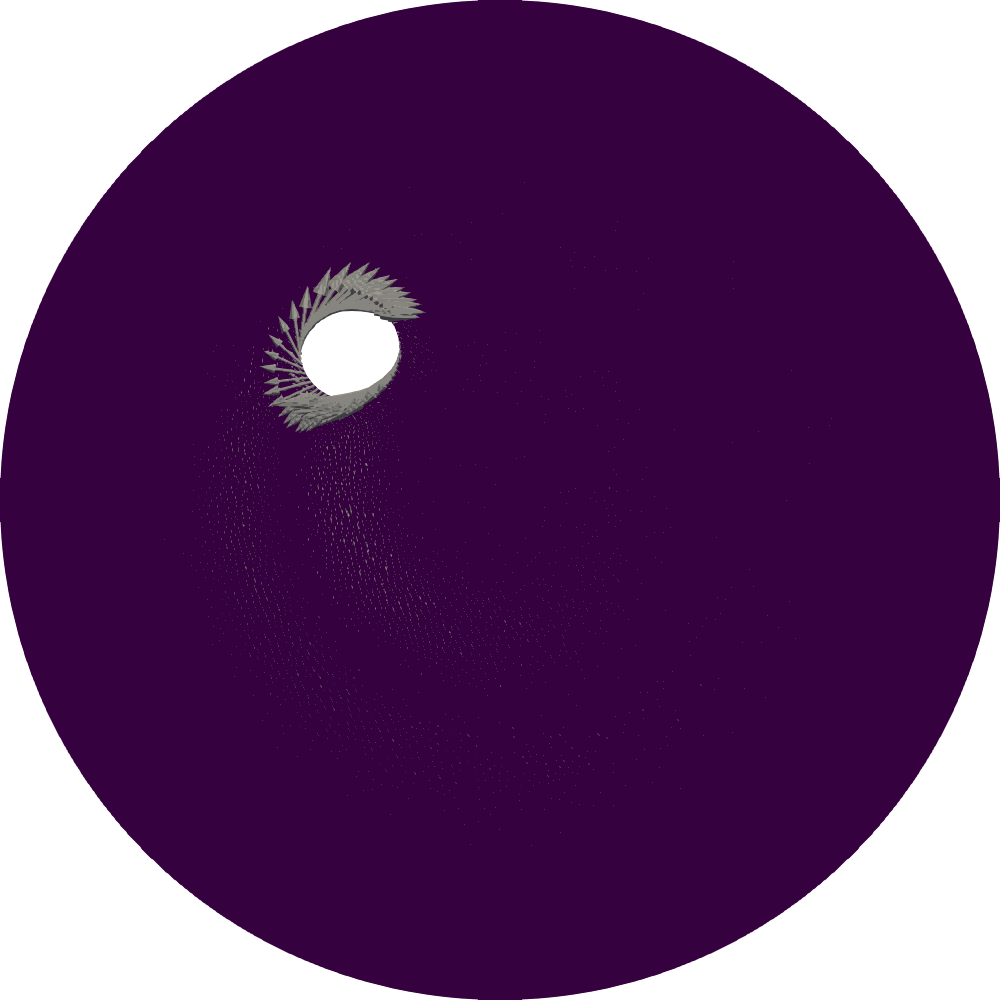}
    \includegraphics[width=0.24\linewidth]{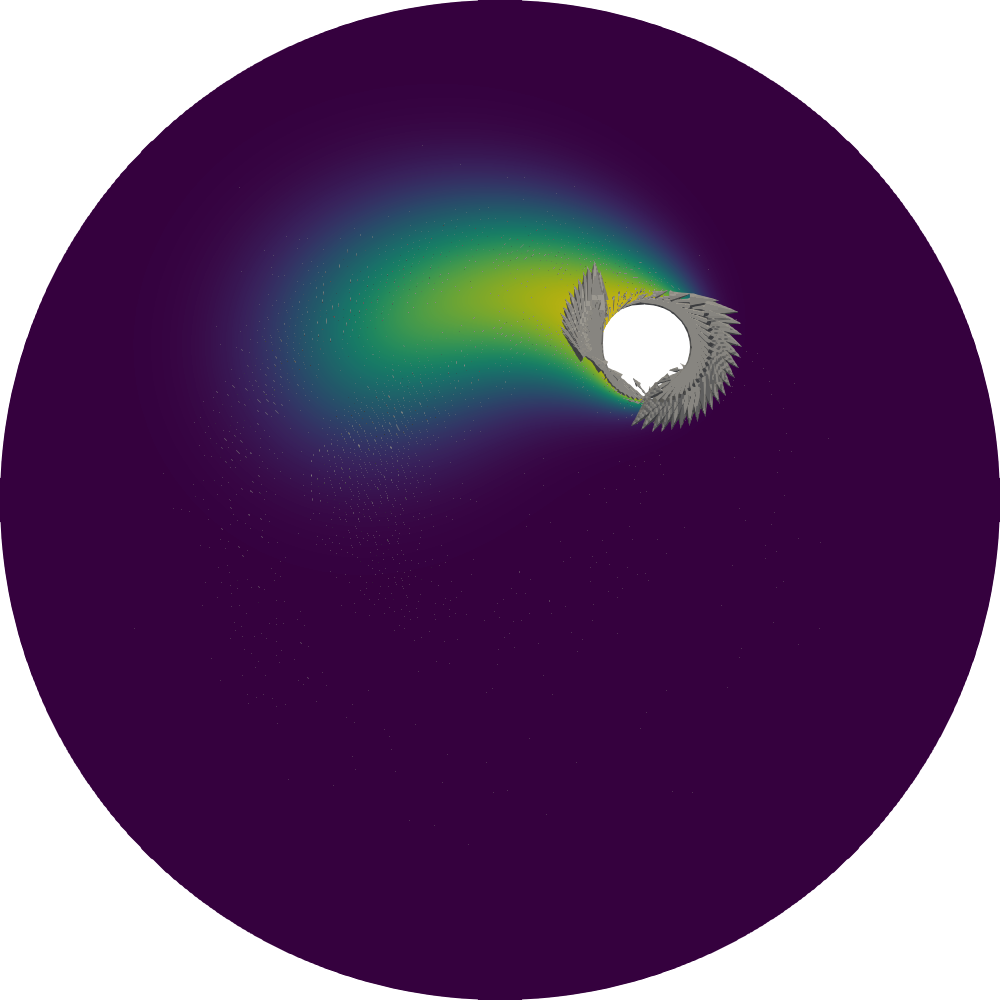}
    \includegraphics[width=0.24\linewidth]{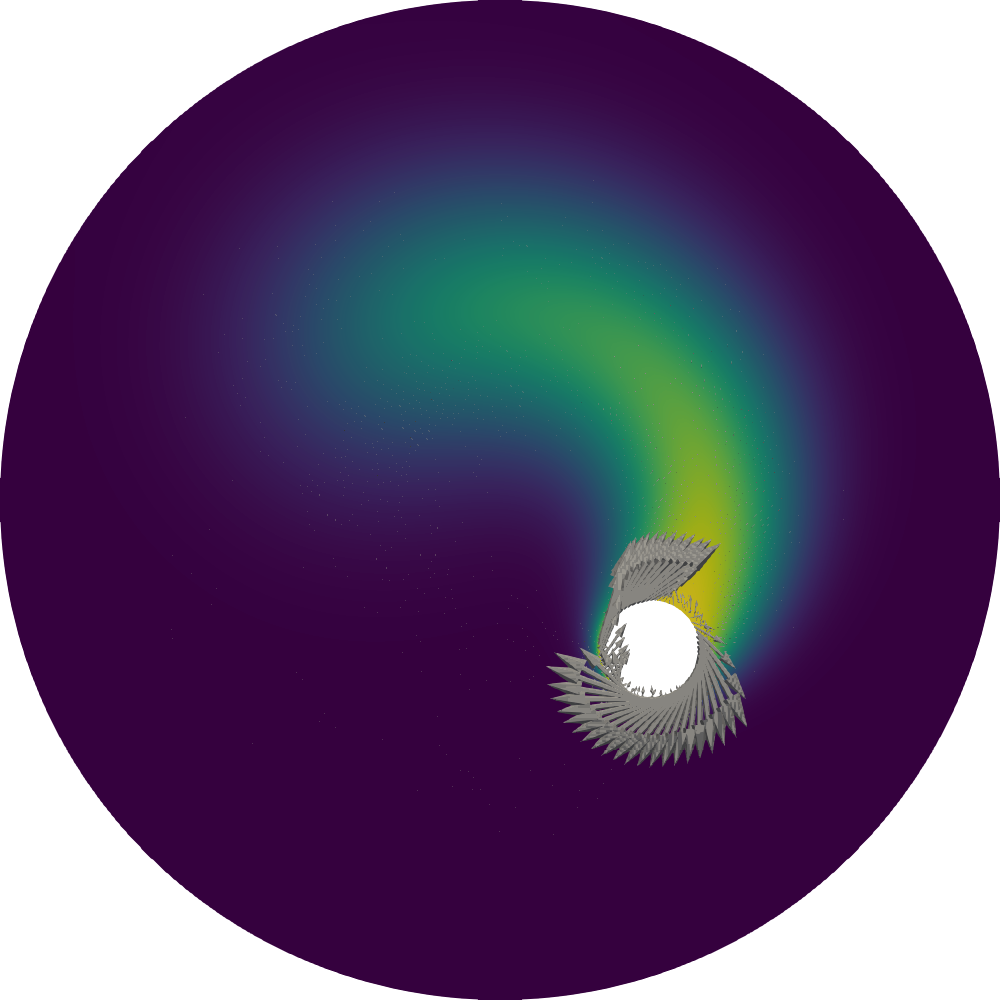}
    \includegraphics[width=0.24\linewidth]{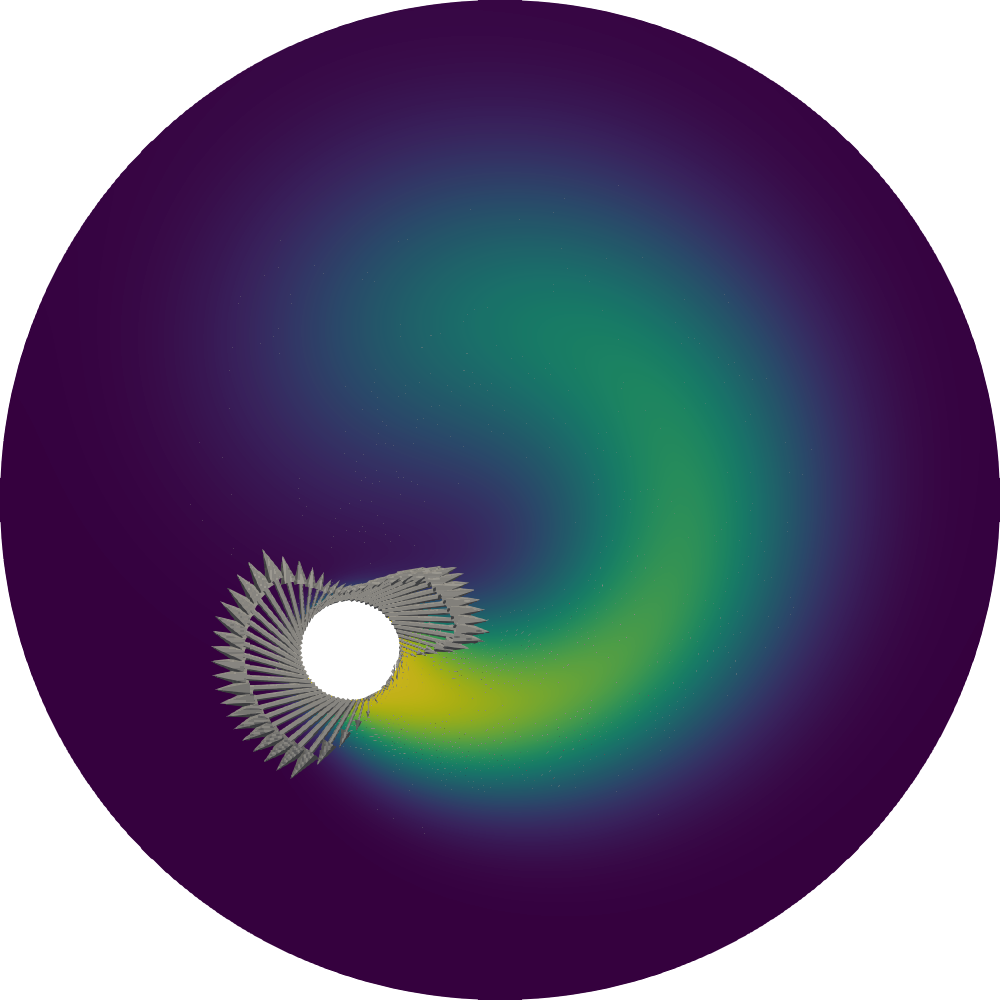}
    \caption{Evolution of the state variable and shape derivative for the time-dependent shape problem. The shape gradients are scaled with $0.1$.}\label{fig:timedep}
  \end{subfigure}
    \begin{subfigure}{\linewidth}
    \includegraphics[width=0.24\linewidth]{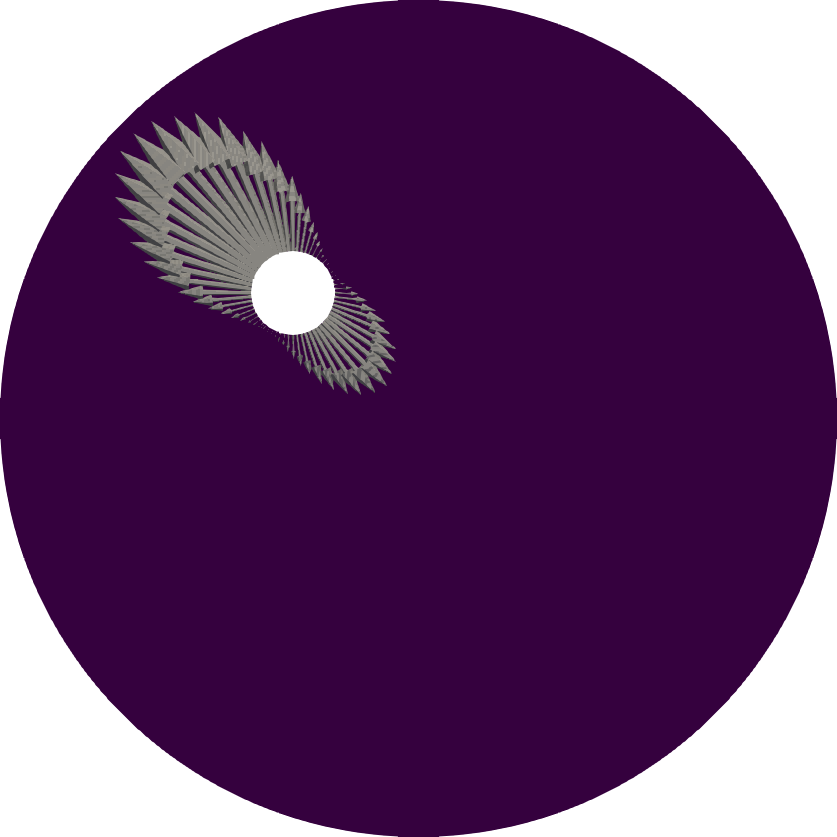}
    \includegraphics[width=0.24\linewidth]{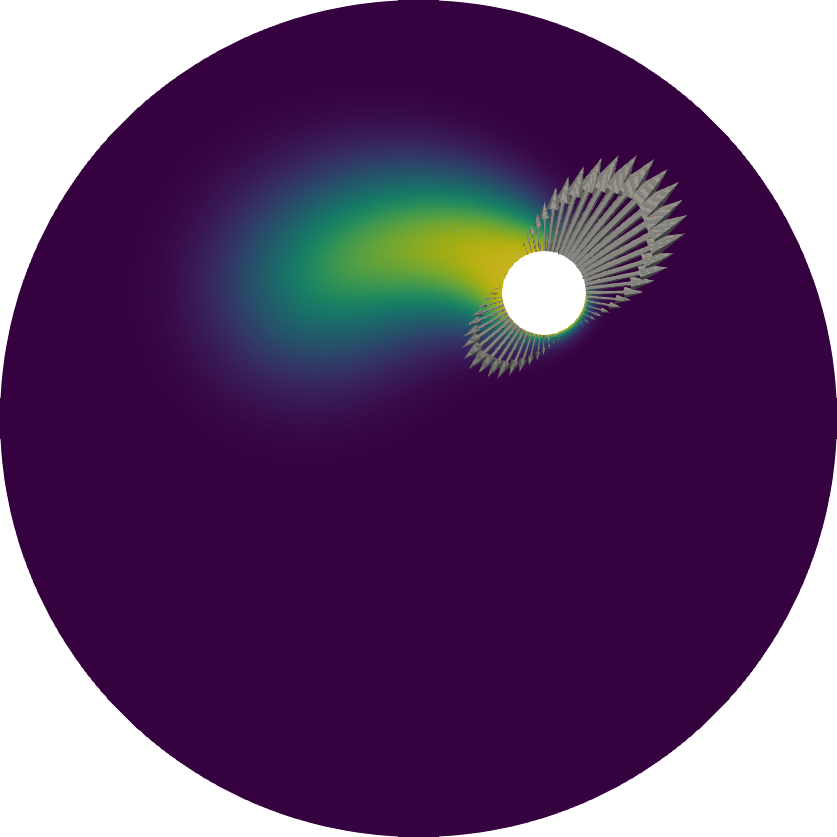}
    \includegraphics[width=0.24\linewidth]{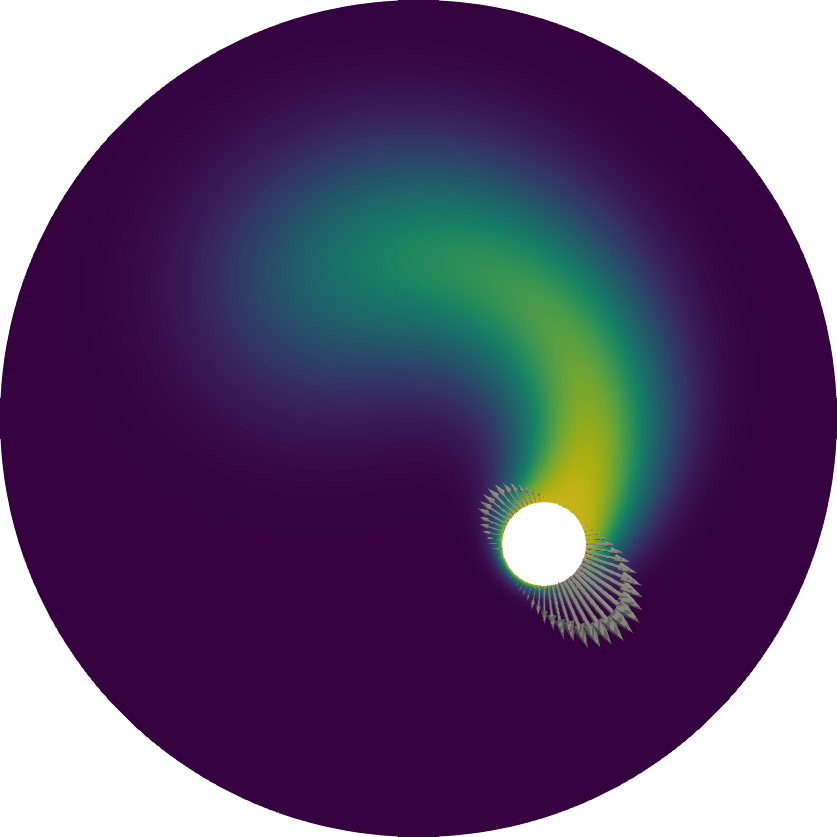}
    \includegraphics[width=0.24\linewidth]{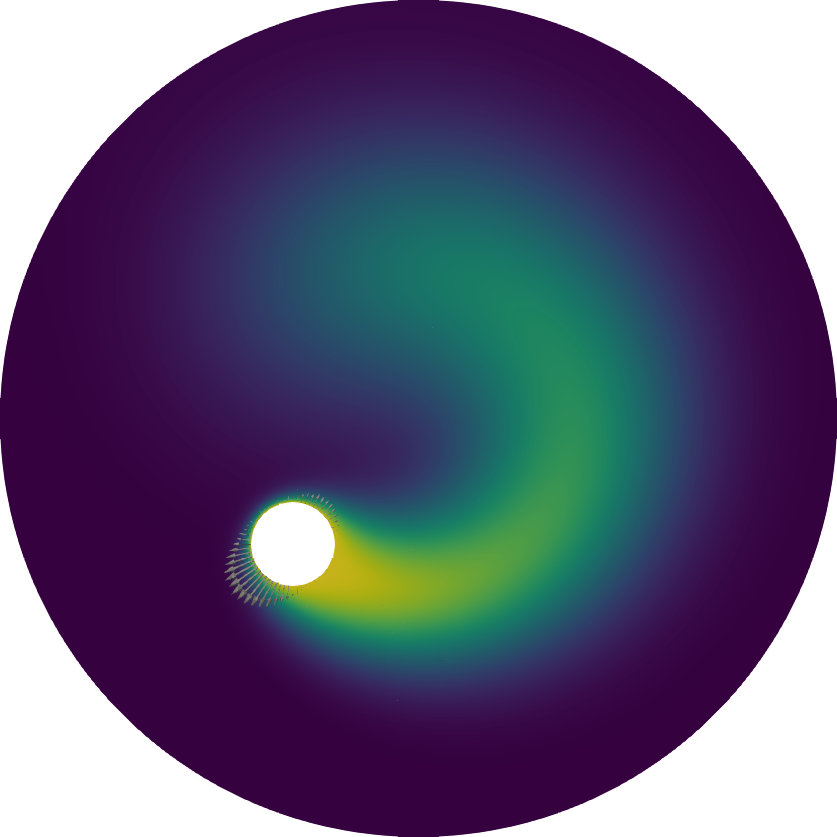}
      \caption{Evolution of the state and shape derivative for the problem with the domain movement decomposed into a static and variable component. The shape gradient is scaled with $0.025$.}\label{fig:base}
  \end{subfigure}
  \caption{Evolution of the state variable and shape derivatives for $T=0,1,2,3$.}\label{fig:tubederivatives}
\end{figure}

The shape gradients for the two approaches \cref{code:approach:1} (left) and \cref{code:approach:2} (right) is visualized in \cref{fig:timedep} and \cref{fig:base} for $t = 0,1,2,3$ using the $l2$ Riesz representation of the gradient.
The key difference is that in the first approach (\cref{code:approach:1}), the rotation of the obstacle is not differentiated through in the shape derivative, and the gradient direction is not the direction of rotation. For the second approach (\cref{code:approach:2}), the differentiation algorithm respects that the obstacle always rotates with a given speed, and the gradient is therefore in the direction of the outer normal, making the heating obstacle wider, emitting more heat.

To verify the algorithmic differentiation algorithm, one can perform Taylor-tests of the reduced functional $\hat J$. This test is based on the fact that
\begin{linenomath}
\begin{subequations}\label{eq:residuals}
\begin{align}
  R_0=&\vert \hat J(\Theta + h \delta \Theta) - \hat J(\Theta)\vert \rightarrow 0 && \text{ at } \mathcal{O}(h),\\
  R_1=&\vert \hat J(\hat\Theta + h \delta \Theta) - \hat J(\Theta)-h\nabla \hat J \cdot (\delta\Theta)\vert\rightarrow 0 && \text{ at } \mathcal{O}(h^2)\\
  R_2=&\vert \hat J(\Theta + h \delta \Theta) - \hat J(\Theta)-h\nabla \hat J \cdot (\delta\Theta)
  - h^2(\delta\Theta)^T\cdot \nabla^2 \hat J \cdot (\delta\Theta)\vert\rightarrow 0 && \text{ at } \mathcal{O}(h^3),
\end{align}
\end{subequations}
\end{linenomath}
where $\Theta=(\hat\theta_0,\cdots,\hat\theta_N)$, $\partial\Theta=(\partial\hat\theta_0,\cdots,\partial\hat\theta_N)$.
This is done in dolfin-adjoint by calling
\begin{python}
  convergence_summary = taylor_to_dict(Jhat, thetas, dthetas)
\end{python}
where \pythoninline{dthetas} is a list of $N+1$ perturbation vectors used in the Taylor test.
Choosing the test directions $\delta\hat\theta_i=(1-x^2-y^2 \quad 1-x^2-y^2)$ yields \cref{tab:tubes} and \cref{tab:tubes:base} for the two different approaches of choosing the control variable. The computational domain consists of $7854$ cells.
\begin{table}[!ht]
  \centering
 \begin{tabular}{ c  c  c  c  c  c  c }
 \hline
  & $R_0$ & \textbf{Rate} & $R_1$ & \textbf{Rate} & $R_2$ & \textbf{Rate} \\ \hline
 $h$   & $3.55\cdot 10^{-3}$ & $-$    & $3.79\cdot 10^{-4}$ & $-$
 & $2.37\cdot 10^{-7}$ & $-$
 \\ \hline
 $h/2$ & $1.87\cdot 10^{-3}$ &$\mathbf{0.93}$& $9.47\cdot 10^{-5}$ &$\mathbf{2.00}$
 & $2.96\cdot 10^{-8}$ &$\mathbf{3.00}$
 \\  \hline
 $h/4$ & $9.56\cdot 10^{-4}$ &$\mathbf{0.96}$ & $2.37\cdot 10^{-5}$ &$\mathbf{2.00}$
 & $3.71\cdot 10^{-9}$ &$\mathbf{3.00}$
 \\  \hline
 $h/8$ & $4.85\cdot 10^{-4}$ &$\mathbf{0.98}$ & $5.92\cdot 10^{-6}$ &$\mathbf{2.00}$
 & $4.64\cdot 10^{-10}$ &$\mathbf{3.00}$
 \\  \hline
 \end{tabular}
 \caption{Residuals and convergence rates for forward problem where the movement only consists of rotation, which is the control variable.}\label{tab:tubes}
\end{table}
\begin{table}[!ht]
  \centering
 \begin{tabular}{ c  c  c  c  c  c  c }
 \hline
  & $R_0$ & \textbf{Rate} & $R_1$ & \textbf{Rate} & $R_2$ & \textbf{Rate} \\ \hline
 $h$   & $6.04\cdot 10^{-2}$ & $-$    & $6.46\cdot 10^{-5}$ & $-$
 & $5.52\cdot 10^{-8}$ & $-$
 \\ \hline
 $h/2$ & $3.02\cdot 10^{-2}$ &$\mathbf{1.00}$ & $1.62\cdot 10^{-5}$ &$\mathbf{2.00}$
 & $6.90\cdot 10^{-9}$ &$\mathbf{3.00}$
 \\  \hline
 $h/4$ & $1.51\cdot 10^{-2}$ &$\mathbf{1.00}$ & $4.04\cdot 10^{-6}$ &$\mathbf{2.00}$
 & $8.63\cdot 10^{-10}$ &$\mathbf{3.00}$
 \\  \hline
 $h/8$ & $7.56\cdot 10^{-3}$ &$\mathbf{1.00}$ & $1.01\cdot 10^{-6}$ &$\mathbf{2.00}$
 & $1.08\cdot 10^{-10}$ &$\mathbf{3.00}$
 \\  \hline
 \end{tabular}
  \caption{Residuals and convergence rates for forward problem where the movement is decomposed in to a fixed motion (rotation), and the motion that is the control variable.}\label{tab:tubes:base}
\end{table}
Finally, we consider the performance of the automatically computed derivatives, and the corresponding adjoint equations.
A comparison of the run-time for the forward, backward and second order adjoint equations are shown in \cref{tab:performance,tab:performance:base} for the two different setups of the problem. In addition to these timings, we compared the run-time of the forward problem with and without the overloading actions in dolfin-adjoint. The overloading actions increased the forward run-time with less than $5$ percent.
\begin{table}[!ht]
  \centering
\begin{tabular}{lcc}
  \hline Operation & Run-time $(s)$  & \textbf{Rate}  \\ \hline
  Forward problem &$116.64$ & -  \\ \hline
  First order derivative (Adjoint problem) & $87.96$ & $\mathbf{0.75}$\\ \hline
  Second order derivative (TLM \& 2nd adjoint problem)& $554.40$ & $\mathbf{4.75}$\\ \hline
\end{tabular}
  \caption{Computational time for different operations for the case where the rotational rotational motion is not differentiated through, and the PDE $F_S$ is not annotated. We used a computational domain consisting of $30 886$ elements, and an end time $T=4$.}\label{tab:performance}
\end{table}
\begin{table}[!ht]
  \centering
\begin{tabular}{lcc}
  \hline Operation & Run-time $(s)$  & \textbf{Rate}  \\ \hline
  Forward problem &$192.90$ & -  \\ \hline
  First order derivative (Adjoint problem) & $178.58$ & $\mathbf{0.93}$\\ \hline
  Second order derivative (TLM \& 2nd adjoint problem)&$754.82$ & $\mathbf{3.91}$\\ \hline
\end{tabular}
\caption{Computational time for the different operations in dolfin-adjoint, when the movement is decomposed into two components,
  a fixed movement, and the movement that is the control variable.
We used a computational domain consisting of $30 886$ elements, and an an end time $T=4$.}\label{tab:performance:base}
\end{table}
\section{Numerical Examples}\label{sec:examples}
In this section, we will present two examples, highlighting the new features of dolfin-adjoint.
First, we solve a shape optimization problem for a stationary PDE with an analytic solution. In this example, we investigate different ways of computing the shape gradient with different mesh deformation techniques.

Then, in the second example we illustrate that dolfin-adjoint can compute shape sensitivities of time-dependent and non-linear PDEs with very little overhead to the forward code.
In this example, we consider a functional consisting of the drag and lift coefficients of an obstacle subject to a Navier-Stokes fluid flow.

\subsection{Pironneau benchmark}\label{sec:stokes}
The first example will illustrate how dolfin-adjoint can be used to solve shape optimization problems with a wide range of approaches. We present how to use Riesz representations with appropriate inner product spaces, how to use custom mesh deformation schemes, as well as how to differentiate through the mesh deformation scheme.

We consider the problem of minimizing the dissipated fluid energy in a channel with a solid obstacle, with the governing equations being the Stokes equations.
This problem has a known analytical solution~\cite{pironneau1974optimum}, an object shaped as an American football with a 90 degree back and front wedge.

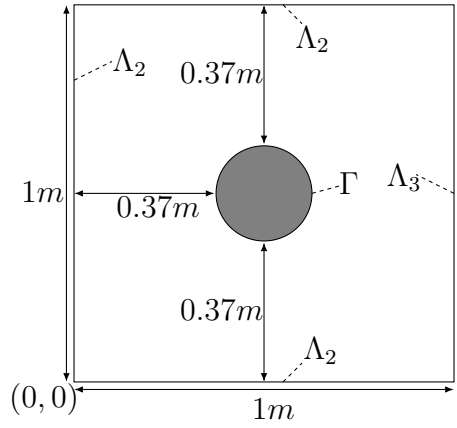
\begin{figure}[!ht]
  \centering
  \begin{overpic}[scale=0.5, grid=false]{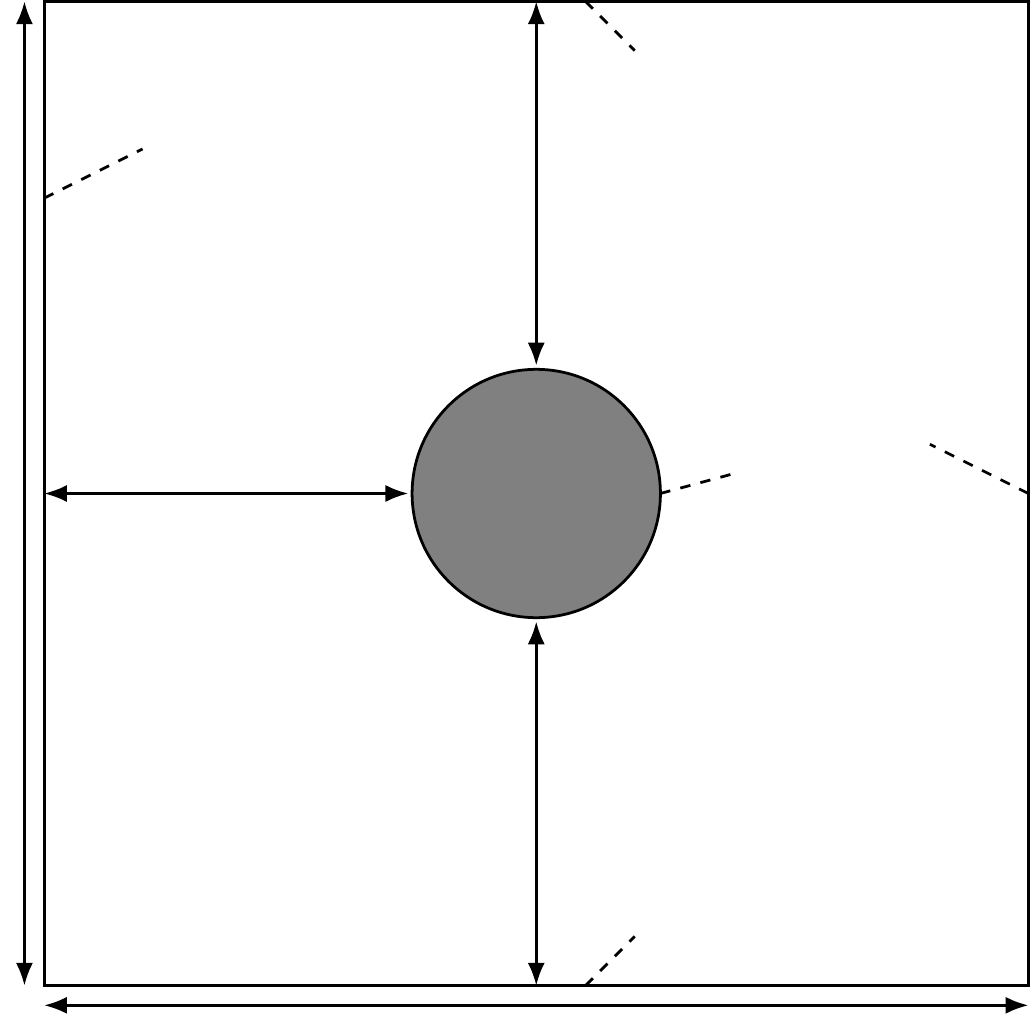}
    \put(-9,50) {$1m$}
    \put(14,83) {$\Lambda_2$}
    \put(60,89) {$\Lambda_2$}
    \put(-12,-2) {$(0,0)$}
    \put(49,-5) {$1m$}
    \put(31,20) {$0.37m$}
    \put(15,46) {$0.37m$}
    \put(31,80) {$0.37m$}
    \put(62,10) {$\Lambda_2$}
    \put(71,52) {$\Gamma$}
    \put(83,54) {$\Lambda_3$}
  \end{overpic}
  \vspace{0.4cm}
	\caption{The setup of the fluid domain for the Pironneau problem, a channel with a circular obstacle in the center.}\label{fig:piro}
\end{figure}

In order to avoid trivial solutions of the optimization problem,
volume and barycenter constraints are added as quadratic penalty terms to the functional.
The optimization problem is written as:
\begin{linenomath}
\begin{align}\label{eq:stokes_opt}
  &\min_{u,p,s} \Int{\Omega(s)}{} \sum_{i,j=1}^2\left(\der{u_i}{x_j}\right)^2 \md x
  + \alpha\left(\mathrm{Vol}(\Omega(s))-\mathrm{Vol}(\Omega) \right)^2
  + \beta \sum_{i=1}^2 \left(\mathrm{Bc}_i(\Omega(s))-\mathrm{Bc}_i(\Omega)\right)^2
\end{align}
\end{linenomath}
subject to
\begin{linenomath}
  \begin{subequations}\label{eq:stokes_state}
\begin{align}
  -\Delta u + \nabla p &= 0 && \text{in } \Omega(s),\\
  \mdiv{u} &= 0 && \text{in } \Omega(s),\\
  u &= 0 && \text{on } \Gamma(s),\\
  u &= g && \text{on } \Lambda_2,\\
  \der{u}{n} + pn &= 0 &&\text{on } \Lambda_3,
\end{align}
\end{subequations}
\end{linenomath}
where $\Omega$ is the unperturbed domain, $\Omega(s)$ the perturbed domain,$\mathrm{Vol}(\Omega) = 1-\Int{\Omega}{}1\md x$ is the volume of the obstacle, $\mathrm{Bc}_i(\Omega) = \left(0.5-\Int{\Omega}{}x_i\md x\right)/\mathrm{Vol}(\Omega)$ is the $i$-th component of the barycenter of the obstacle. The fluid velocity and pressure is denoted $u$ and $p$, respectively.
$\alpha$ and $\beta$ are penalty parameters for the quadratic volume and barycenter penalization.
The unperturbed domain is visualized in \cref{fig:piro}.

The forward problem \eqref{eq:stokes_state} can be implemented in FEniCS as shown in \cref{code:stokes:forward}.
  \begin{python}[caption={Code for solving the Stokes problem \eqref{eq:stokes_state} as a function of a mesh perturbation field $s$. For brevity, we skip the volume and barycenter constraints.},captionpos=b, label={code:stokes:forward},numbers=left]
def forward(mesh, mf, s):
    # Compute initial volume and barycenter
    # ...

    # Move mesh with perturbation s
    ALE.move(mesh, s)

    # Define variational formulation of Stokes problem
    Ve = VectorElement("CG", mesh.ufl_cell(), 2)
    Qe = FiniteElement("CG", mesh.ufl_cell(), 1))
    VQ = FunctionSpace(mesh, Ve*Qe)
    (u, p) = TrialFunctions(VQ)
    (v, q) = TestFunctions(VQ)
    a = inner(grad(u), grad(v))*dx - div(u)*q*dx - div(v)*p*dx
    l = inner(Constant((0,0)), v)*dx

    # Create boundary conditions
    markers = {"inflow":1, "outflow":2, "walls":3, "obstacle": 4}
    g = {"inflow": Expression(("sin(pi*x[1])","0"), element=Ve),
           "walls": (0,0), "obstacle": (0,0)}
    bcs = [DirchletBC(V.sub(0), g[key], mf, markers[key])
           for key in bcs.keys()]

    # Solve variational problem
    w = Function(VQ, name="Mixed State Solution")
    solve(a==l, w, bcs=bcs)
    u, p = w.split()

    J = assemble(inner(grad(u), grad(u))*dx)
    # Add barycenter and volume constraints
    # ...
    return J
  \end{python}
The shape derivative can then be obtained with the additions shown in \cref{code:stokes:overhead}.
\begin{python}[caption={The necessary additions dolfin-adjoint requires to compute shape sensitivities of the functional \eqref{eq:stokes_opt}.},captionpos=b, label={code:stokes:overhead}]
from dolfin import *
from dolfin_adjoint import *
mesh = Mesh("file.xml")
mf = MeshFunction("size_t", mesh, "facet_function.xml")
S = VectorFunctionSpace(mesh, "CG", 1)
s = Function(S)
J = forward(mesh, mf, s)
Jhat = ReducedFunctional(J, Control(s))
dJds = Jhat.derivative()
  \end{python}

\subsubsection{Custom mesh deformation schemes}
As for the documented example in \cref{sec:FEniCS}, the shape derivative will have its main contributions on the boundary. To use a $L^2(\Omega)$ Riesz representation of the shape derivative to perturb the domain will lead to a degenerated mesh. Similarly, the $H^1(\Omega)$ Riesz representation often yield degenerate meshes for large deformations.

Therefore, we introduce a mesh deformation scheme that consist of bilinear forms $a(d,s)=\mathrm{d}J(\Omega(\hat\theta))[s]$. There exists a wide variety of such schemes, for instance a linear elasticity equation with spatially varying Lamé parameters~\cite{schulz2016computational}, restricted mesh deformations~\cite{etling2018first} and convection-diffusion equations using Eikonal equations for distance measurements~\cite{schmidt2014two}.

To solve the optimization problem \eqref{eq:stokes_opt}, we use Moola~\cite{moola}. What distinguishes Moola from many optimization packages, is that it uses the native inner products to determine search directions and convergence criteria. This means that for functions living in $H(\Omega)$, it uses the inner product $(u,v)_{\Omega}=\int_{\Omega}u \cdot v \md x$.
\cref{code:moola} illustrates how to use the Moola interface in combination with \cref{code:stokes:overhead}, using a Newton-CG solver with a custom Riesz representation.
  \begin{python}[caption={Code illustrating how to use Moola in combination with dolfin-adjoint. To use custom Riesz representations in Moola, the DolfinPrimalVector has do be supplied with a map from the primal to dual space and its inverse.},captionpos=b, label={code:moola}]
import moola
problem = MoolaOptimizationProblem(Jhat)
class CustomRieszMap(object):
def __init__(self):
    # ...
    def primal_map(self, x, b):
        # ...
    def dual_map(self, x):
        # ...
f_moola = moola.DolfinPrimalVector(s, riesz_map=CustomRieszMap())
solver = moola.NewtonCG(problem, f_moola)
s_opt = solver.solve()
  \end{python}

\subsubsection{Custom descent schemes}
Using Riesz representations without mesh-quality control can lead to inverted/degenerated elements, and the tolerances for the Newton-CG method has to be chosen appropriately. An alternative approach would be to use a restricted gradient descent scheme, as presented by~\cite{etling2018first}, where the Riesz representation has additional restrictions, as well as a mesh quality control check in the descent scheme.

\subsubsection{Differentiation of the deformation schemes}\label{sec:meshdef}
In the two previous approaches, the shape gradient is first computed, then corresponding mesh deformation is computed through mesh deformation schemes, with or without restrictions.
A third approach is to differentiate through the mesh deformation scheme.
This is illustrated by employing a slightly modified version of the elasticity equation presented in \cite{schulz2016computational}.
We rewrite the optimization problem \eqref{eq:stokes_opt} as
\begin{linenomath}
\begin{subequations}
\begin{align}\label{eq:stokes_opt2}
  &\min_{u,p,s,h} \Int{\Omega(s(h))}{} \sum_{i,j=1}^2\left(\der{u_i}{x_j}\right)^2 \md x
  + \alpha\left(\mathrm{Vol}(\Omega(s(h)))-\mathrm{Vol}(\Omega) \right)^2
  + \beta \sum_{i=1}^2 \left(\mathrm{Bc}_i(\Omega(s(h)))-\mathrm{Bc}_i(\Omega)\right)^2
\end{align}
\end{subequations}
\end{linenomath}

subject to \cref{eq:stokes_state} and
\begin{linenomath}
\begin{subequations}\label{eq:mesh-deformation}
\begin{align}
  \mdiv{\sigma}&=0\quad\text{in } \Omega_0,\\
  s &= 0\quad\text{on } \Lambda_1\cup\Lambda_2\cup\Lambda_3,\\
  \der{s}{n} &= h \quad \text{on } \Gamma,
\end{align}
\end{subequations}
\end{linenomath}
where the stress tensor $\sigma$ and strain tensor $\epsilon$ is defined as
\begin{linenomath}
\begin{subequations}
\begin{align}
  \sigma &:= \lambda_{elas}\mathrm{Tr}(\epsilon) + 2\mu_{elas}\epsilon,\\
  \epsilon &:= \half(\nabla s + \nabla s^T).
\end{align}
\end{subequations}
\end{linenomath}
As in~\cite{schulz2016computational}, we set the Lamé parameters as $\lambda_{elas}=0$ and let $\mu_{elas}$ solve
\begin{linenomath}
\begin{subequations}\label{eq:lame}
\begin{align}
  \Delta \mu_{elas} &= 0 \quad \text{in } \Omega_0,\\
  \mu_{elas} &= 1 \quad \text{on } \Lambda_1\cup\Lambda_2\cup\Lambda_3,\\
  \mu_{elas} &= 500 \quad \text{on } \Gamma.
\end{align}
\end{subequations}
\end{linenomath}
This approach can be though of as finding the boundary stresses that deforms the mesh such the energy dissipation in the fluid is minimized.
Using this approach, a Riesz representation of the control $h\in L^2(\partial\Omega)$ suffices, as the mesh deformation is a function of the control variable.

This is implemented in dolfin-adjoint using the \pythoninline{BoundaryMesh} class and the \pythoninline{transfer_from_boundary} function. An outline of the implementation is given in \cref{code:stokes:overhead2}.
\begin{python}[caption={Overhead for differentiating through mesh deformations, where the design variable is a finite element function defined only at the boundary of the computational domain.},captionpos=b, label={code:stokes:overhead2}]
from dolfin import *
from dolfin_adjoint import *
# Load mesh and mesh function from file
# ...
# Define the BoundaryMesh and the design variable
b_mesh = BoundaryMesh(mesh, "exterior")
S_b = VectorFunctionSpace(b_mesh, "CG", 1)
h = Function(S_b, name="Design")

# Transfer values from the FunctionSpace on the BoundaryMesh,
# to the FunctionSpace on the Mesh (to be used in the
# variational formulation).
h_V = transfer_from_boundary(h, mesh)

# Deform mesh according to the extension
s = mesh_deformation(mesh, mf, h_V)

# Solve forward problem and define reduced functional
J, u = forward(mesh, mf, s)
Jhat = ReducedFunctional(J, Control(h))
\end{python}

\subsubsection{Results}
\begin{figure}[!ht]
\centering
  \begin{subfigure}{0.4\linewidth}
    \includegraphics[width=\linewidth]{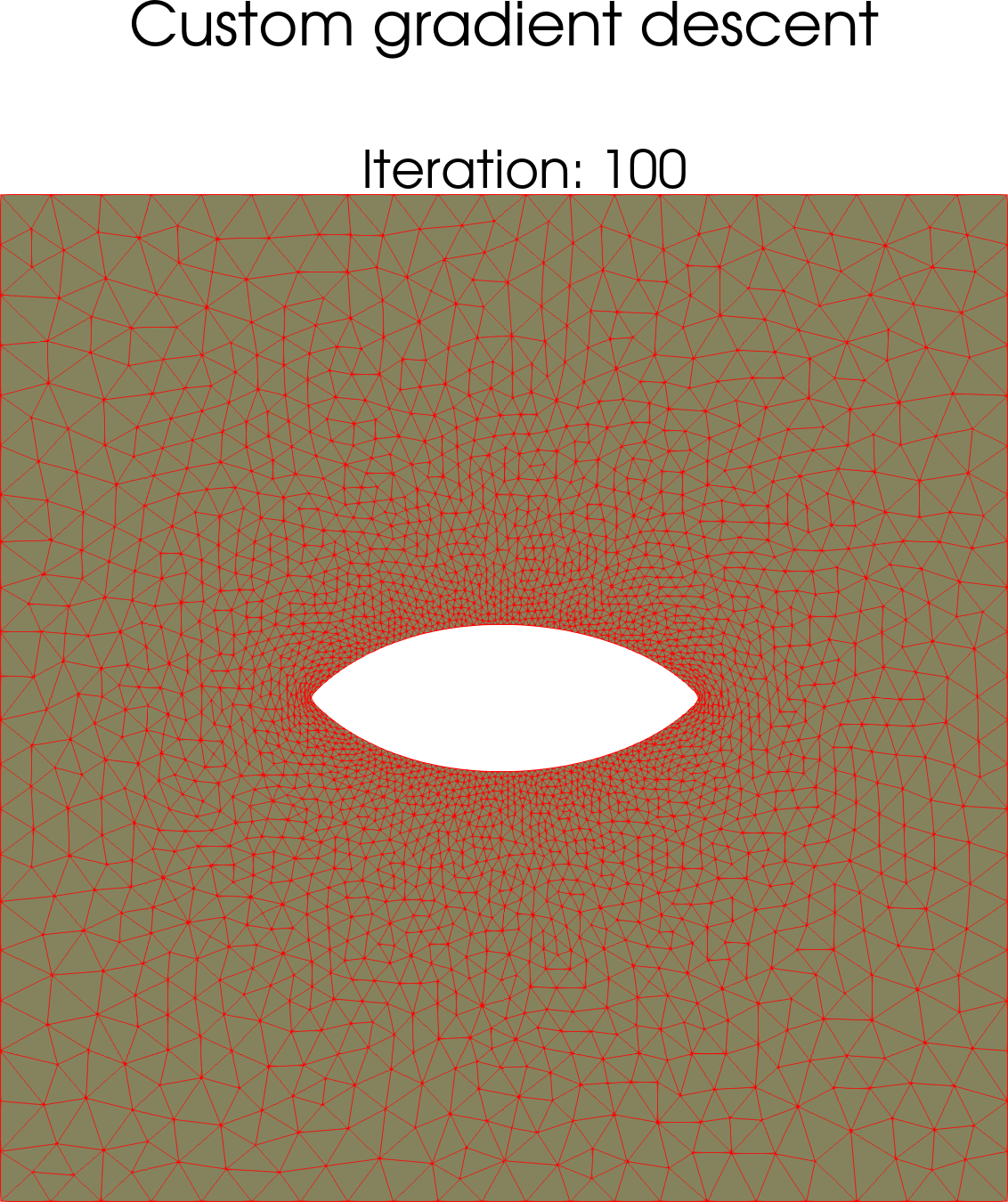}
    \caption{Mesh deformation after 100 iterations of a custom steepest descent algorithm, using a restricted gradient representation.}\label{roland}
  \end{subfigure}
  \begin{subfigure}{0.4\linewidth}
      \includegraphics[width=\linewidth]{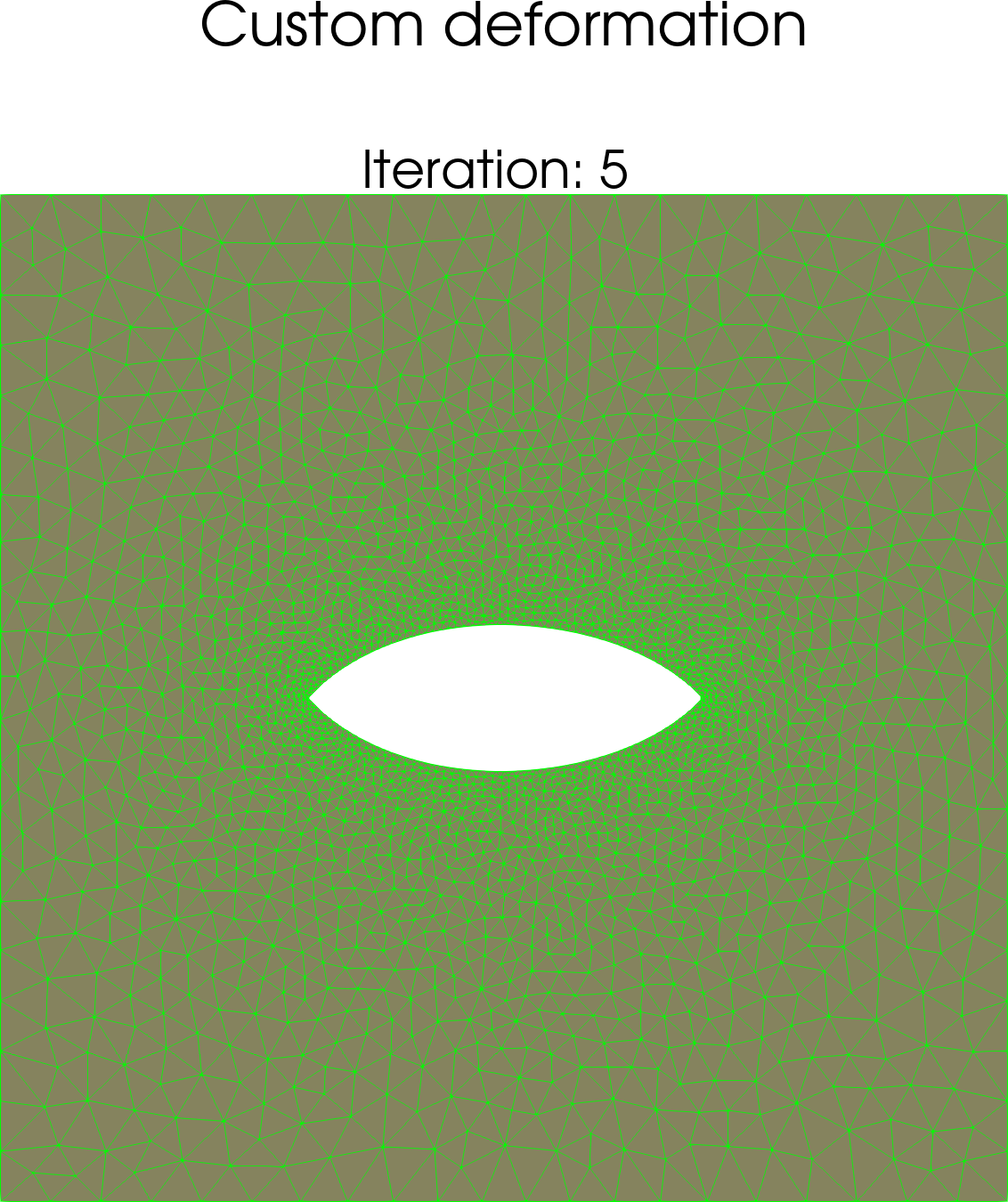}
  \caption{Mesh deformation after 5 iterations of a Newton-CG algorithm using a custom deformations scheme as a Riesz representation in Moola.}\label{volker}
  \end{subfigure}
  \begin{subfigure}{0.4\linewidth}
      \includegraphics[width=\linewidth]{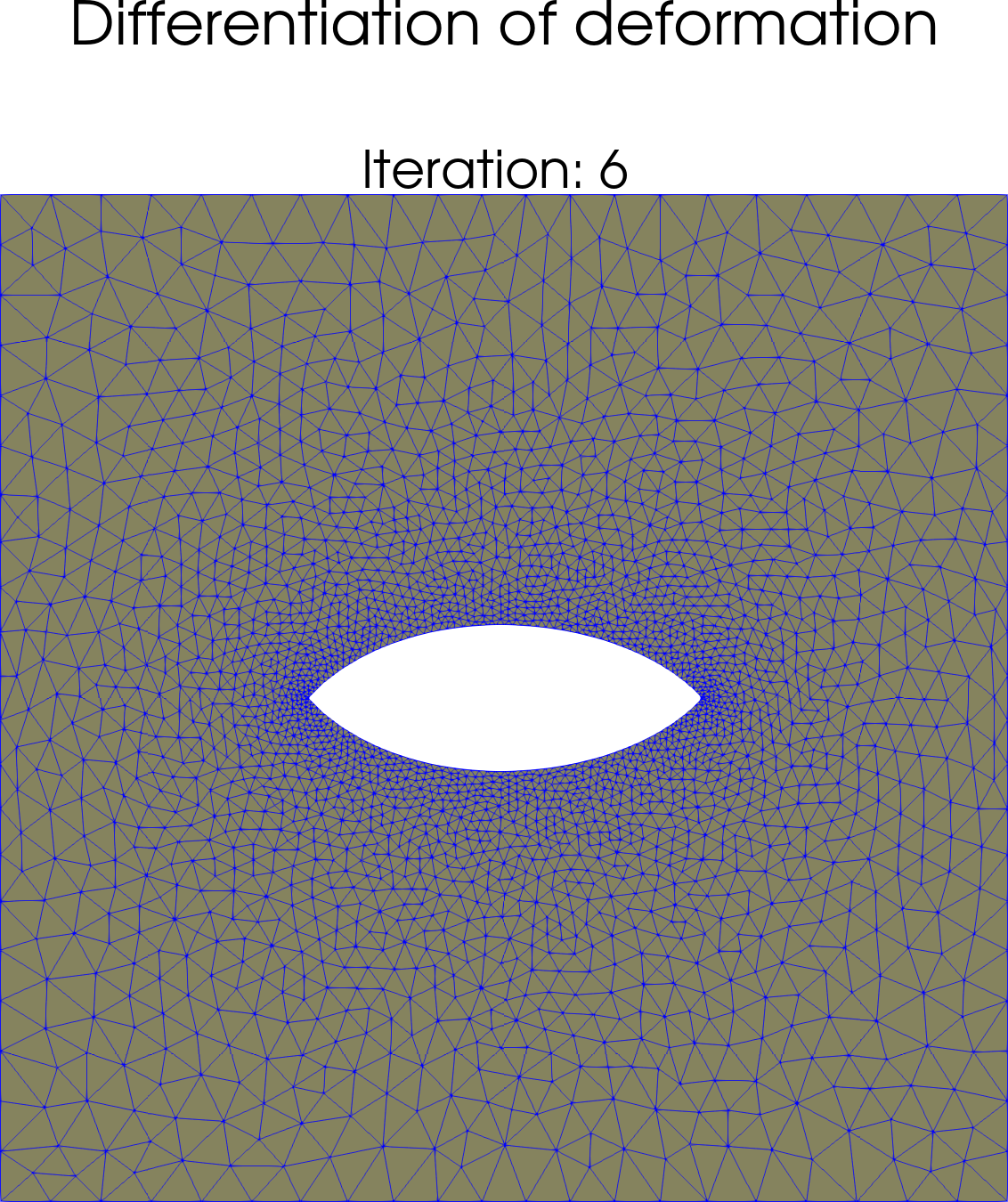}
  \caption{Mesh deformation after 6 iterations of a Newton-CG algorithm differentiating through the mesh deformation scheme.}\label{dokken}
  \end{subfigure}
  \begin{subfigure}{0.4\linewidth}
    \vspace{2.7cm}
      \includegraphics[width=\linewidth]{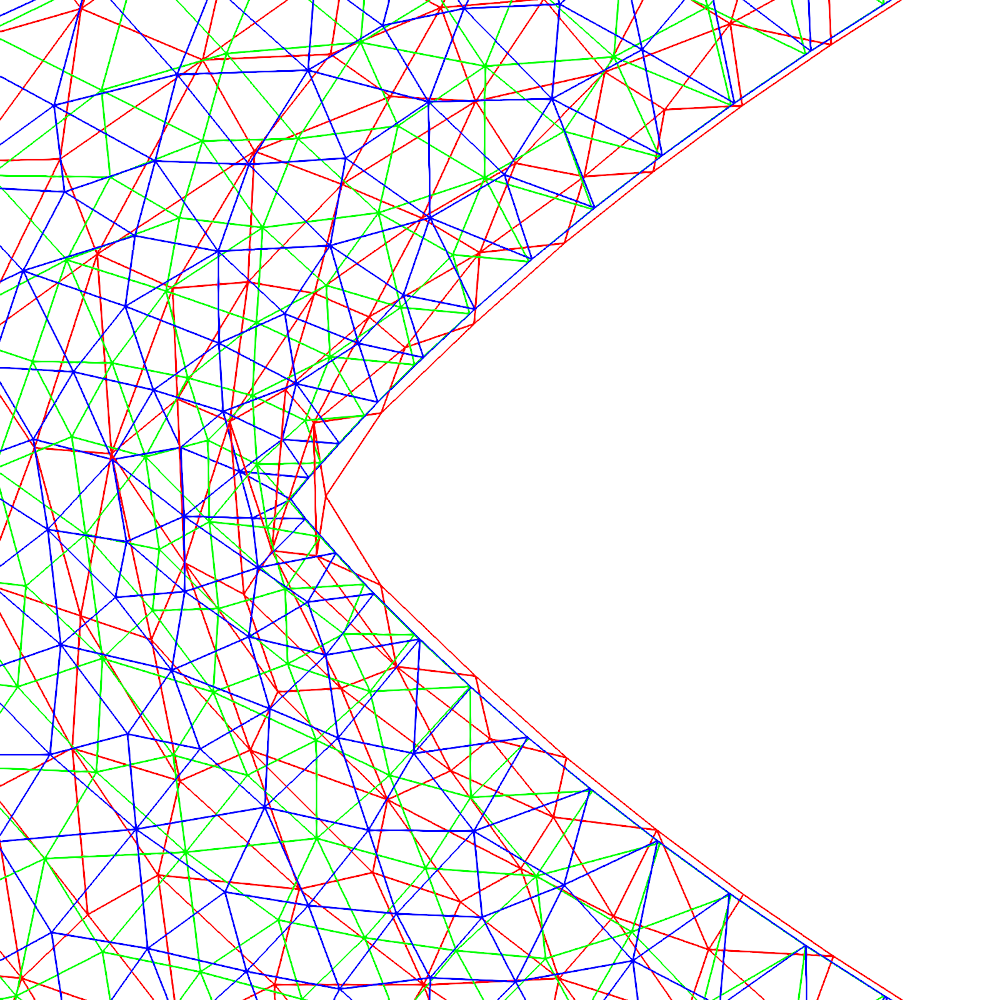}
  \caption{Comparison of the front wedge of the meshes for each deformation scheme at their final iteration. The two Newton-CG based algorithms have almost the exact boundary representation, but with different interior node displacement.
  }\label{compare}
  \end{subfigure}
  \caption{Comparison of the different methods of solving the shape optimization, using (\protect\subref{roland}) A custom gradient representation combined with a Newton-CG solver, (\protect\subref{volker}) a custom steepest descent algorithm with mesh quality checks and restrictions on the gradient representation,  (\protect\subref{dokken}) differentiation through the mesh deformation scheme, only using the boundary nodes as design parameters.}\label{fig:NCG:compare}
\end{figure}

In \cref{fig:NCG:compare}, we compare the three approaches described in the last three sections.
The custom steepest descent algorithm was manually terminated after 100 iterations.
The custom deformation scheme was terminated when the $L^2(\Omega)$ norm of the gradient representation was less than $9\cdot 10^{-3}$ (with a total of 16 conjugate gradient iterations).
As the gradient representation for the custom deformation scheme is not discretely consistent, a lower termination criteria can not be set.
For the differentiation through the mesh deformation scheme, the gradient termination criteria of $5\cdot10^{-6}$ was reached after 6 iterations (with a total of 71 conjugate gradient iterations).
The drag was reduced from $24.3019$ to  $20.5397$ for the custom gradient descent scheme, $20.5393$ for the custom deformation scheme and $20.5385$ for the differentiation through mesh deformation.

\FloatBarrier

\subsection{Non-linear time-dependent Navier-Stokes equations}\label{sec:dfg-3}
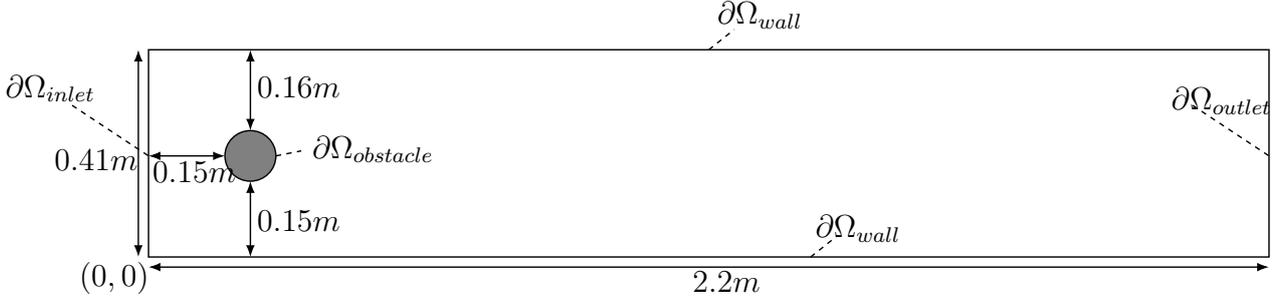
\begin{figure}[!ht]
  \centering
  \begin{overpic}[width=0.95\linewidth,grid=false]{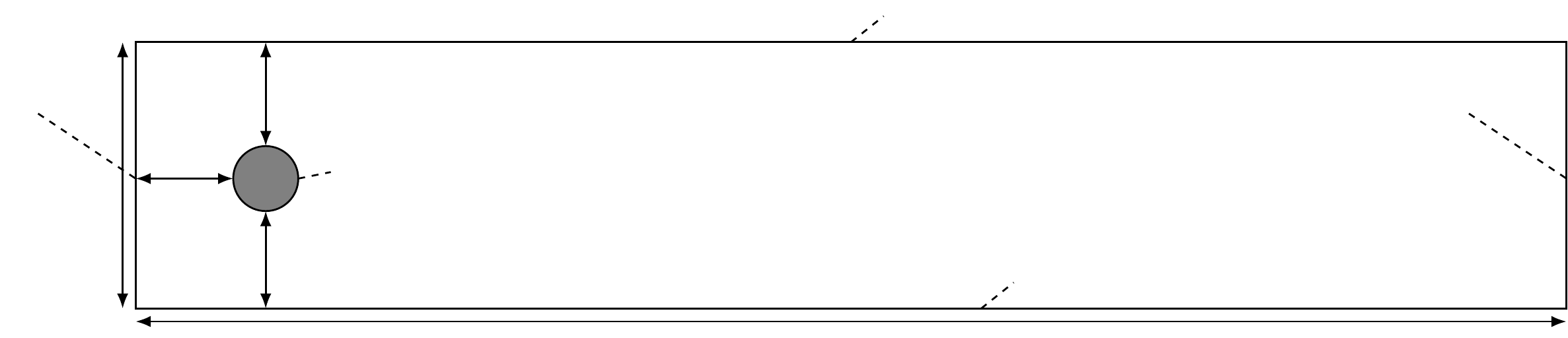}
    \put(3,-0.5) {$(0,0)$}
    \put(-3,15) {$\partial\Omega_{inlet}$}
    \put(1,9) {$0.41m$}
    \put(55,21) {$\partial\Omega_{wall}$}
    \put(53,-1) {$2.2m$}
    \put(17.5,4) {$0.15m$}
    \put(17.5,15) {$0.16m$}
    \put(9,8) {$0.15m$}
    \put(22,10) {$\partial\Omega_{obstacle}$}
    \put(63,3.5) {$\partial\Omega_{wall}$}
    \put(92,14) {$\partial\Omega_{outlet}$}
    \end{overpic}
  \caption{The computational domain for the DFG-3 benchmark described in~\ref{sec:dfg-3}.}\label{fig:NS_domain}
\end{figure}

The aim of this example is to compute shape derivatives for the Featflow DFG-3 benchmark~\cite{FEATFLOW}.
This example is challenging because the Navier-Stokes problem consists of a transient, non-linear equation with a non-trivial coupling between the velocity and pressure field.
We write the Navier-Stokes equations on the following form: Find the velocity $u$ and pressure $p$ such that
\begin{linenomath}
  \begin{subequations}\label{eq:navier-stokes}
  \begin{align}
  \der{u}{t} + u\cdot \nabla u - \nu \Delta u + \nabla p &=0
  &&\text{ in } \Omega\times (0,T],\\
    \nabla \cdot u &= 0 && \text{ in } \Omega\times (0,T],\\
      u(x,t) &= (0,0) &&\text{ on } \partial\Omega_{walls}\cup\partial\Omega_{obstacle}\times (0,T]\\\
        u(x,t)&=\left(\frac{6\sin(\frac{\pi t}{T})x_1(H-x_1)}{H^2},0\right) &&
        \text{ on } \partial\Omega_{inlet}\times (0,T],\\
          \nu \der{u}{n} &= pn && \text{ on } \partial\Omega_{outlet},\\
          u(x,0) &= 0 &&\text{ in } \Omega,
  \end{align}
  \end{subequations}
\end{linenomath}
where $\Omega$ is visualized in \cref{fig:NS_domain}, $\nu=0.001$ the kinematic viscosity, $T=8$ the end time and $H=0.41$ the height of the fluid channel.

For two-dimensional problems, the drag and lift coefficients can be written as~\cite{schafer1996benchmark}
\begin{linenomath}
\begin{subequations}
  \label{eq:cd_cl}
  \begin{linenomath}\begin{align}
          C_D({u},p,t) &= \frac{2}{\rho L U_{mean}^2}\int_{\partial\Omega_{obstacle}}\left(\rho \nu {n} \cdot \nabla u_{t}n_y -pn_x\right)\mathrm{d}s,\label{eq:CD}\\
          C_L({u},p,t) &= -\frac{2}{\rho L U_{mean}^2}\int_{\partial\Omega_{obstacle}}\left(\rho \nu {n} \cdot \nabla u_{t}n_x + pn_y\right)\mathrm{d}s,
    \end{align}
  \end{linenomath}
\end{subequations}
\end{linenomath}
where $n=(n_x,n_y)$ is the outward pointing normal vector, $u_{t}={u}\cdot (-n_y,n_x)$ is the tangential velocity component at the interface of the obstacle $\partial\Omega_{obstacle}$, $U_{mean}=1$ the average inflow velocity, $\rho=1$ the fluid density and $L=0.1$ the characteristic length of the flow configuration.
We chose the functional $J$ as an integrated linear combination of the drag coefficient $C_D$ and the lift coefficient $C_L$.
\begin{linenomath}
\begin{align}\label{eq:J:NS}
  J(u,p) &=  \Int{0}{T}C_D(u,p,t)-C_L(u,p,t)\md t.
\end{align}
\end{linenomath}

As in \cref{sec:meshdef}, we define a perturbation of the computational domain $\Omega(s(h))$, where $s$ is the solution of an elasticity equation
\begin{linenomath}
\begin{subequations}
\begin{align}
  \mdiv{\sigma}&=0\quad\text{in } \Omega,\\
  s &= 0\quad\text{on } \partial\Omega \setminus \partial\Omega_{obstacle},\\
  \der{s}{n} &= h \quad \text{on } \partial\Omega_{obstacle},
\end{align}
\end{subequations}
\end{linenomath}
and the Lamé parameters are set in a similar fashion as in \cref{eq:lame}.
The design parameters of this problem is therefore the stress applied to the mesh vertices at the boundary $\partial\Omega_{obstacle}$.


The Navier-Stokes equations \eqref{eq:navier-stokes} are discretized in time using backward Euler method and a time-step of $\Delta t=1/200$.
For the spatial discretization, we use the Taylor-Hood finite element pair. The mixed velocity pressure function space has 14,808 degrees of freedom.
The non-linear problem at each time-level is solved using the Newton method.

The first and second order shape derivatives of \cref{eq:J:NS} is computed with respect to a change in $h$, and is verified with a Taylor-test in a similar fashion as in \cref{sec:FEniCS}. The results are listed~\cref{fig:taylor:NS} and shows the expected convergence rates.
\begin{table}[!ht]
  \centering
 \begin{tabular}{ c  c  c  c  c  c  c }
 \hline
  & $R_0$ & \textbf{Rate} & $R_1$ & \textbf{Rate} & $R_2$ & \textbf{Rate} \\ \hline
 $h$   & $1.16\cdot 10^{-3}$ & $-$    & $4.32\cdot 10^{-7}$ & $-$
 & $7.15\cdot 10^{-10}$ & $-$
 \\ \hline
 $h/2$ & $7.89\cdot 10^{-4}$ &$\mathbf{1.00}$ & $1.08\cdot 10^{-7}$ &$\mathbf{2.00}$
 & $9.0\cdot 10^{-11}$ &$\mathbf{2.99}$
 \\  \hline
 $h/4$ & $3.94\cdot 10^{-4}$ &$\mathbf{1.00}$ & $2.70\cdot 10^{-8}$ &$\mathbf{2.00}$
 & $1.11\cdot 10^{-11}$ &$\mathbf{3.02}$
 \\  \hline
 $h/8$ & $1.97\cdot 10^{-4}$ &$\mathbf{1.00}$ & $6.75\cdot 10^{-9}$ &$\mathbf{2.00}$
 & $1.27\cdot 10^{-12}$ &$\mathbf{3.12}$
 \\  \hline
 \end{tabular}
 \caption{Taylor test showing the zeroth, first and second order Taylor expansion \eqref{eq:residuals} for the functional \eqref{eq:J:NS}. We observe expected convergence rates for each of the expansions.}\label{fig:taylor:NS}
\end{table}

In \cref{tab:NS}, we time the different operations performed by dolfin-adjoint.
The adjoint computation is faster than the forward computation, as the forward computation is non-linear, and require on average $2$ Newton iterations per time-step.
\begin{table}[!ht]
\centering
\begin{tabular}{lcc}
  \hline Operation & Run-time $(s)$  & \textbf{Rate}  \\ \hline
  Forward problem &$465.34$ & -  \\ \hline
  First order derivative (Adjoint problem) & $312.00$ & $\mathbf{0.67}$\\ \hline
  Second order derivative (TLM \& 2nd adjoint problem)&$807.52$ & $\mathbf{1.74}$\\ \hline
\end{tabular}
\caption{Timings of the operations for computing the forward solution, and the first and second order derivative for the Navier-Stokes problem with a total of 14 808 degrees of freedom for the mixed problem.}\label{tab:NS}
\end{table}

The implementation of the mesh deformation scheme consists of 26 lines of Python code. The forward problem consists of 45 lines of code. The IO for reading in meshes and corresponding markers from XDMF is 7 lines of code.
The additional overhead that has to be added to the code to do automatic shape differentiation of the problem is 7 lines of code.
\section{Concluding remarks}\label{sec:conclusion}
In this paper we have presented a new framework for solving PDE constrained
shape optimization problems for transient domains using high-level algorithmic differentiation on the finite element frameworks FEniCS and Firedrake.
We have demonstrated the flexibility of the implementation, by considering several different aspects of shape optimization, as time-dependent, non-linear problems and time-dependent shapes.
We have verified the implementation by solving a shape optimization problem an analytic solution.
Additionally, the automatically computed first and second order shape derivatives have been verified through Taylor expansions.
In this paper, we have limited the presentation to geometries described by first order geometries.
Therefore, an natural extension to the current software would be to support higher order geometries.

\bibliographystyle{plainurl}
\bibliography{bib}{}
\end{document}